\PassOptionsToPackage{numbers,sort&compress}{natbib}
\documentclass[preprint,12pt]{elsarticle}

\makeatletter
\def\ps@pprintTitle{%
  \let\@oddhead\@empty
  \let\@evenhead\@empty
  \let\@oddfoot\@empty
  \let\@evenfoot\@empty
}
\makeatother

\usepackage[english]{babel}
\usepackage[T1]{fontenc}
\usepackage[utf8]{inputenc}

\usepackage[margin=1.1in]{geometry}
\usepackage{lineno}
\linenumbers

\usepackage{amsmath,amssymb,amsthm,mathtools,bm,stmaryrd,bbm,mathrsfs}

\usepackage{algorithm}
\usepackage{algpseudocode}

\renewcommand*\Call[2]{\textproc{#1}(#2)}

\usepackage{tabularx,array,multirow,booktabs,diagbox}
\usepackage{float}
\usepackage[section]{placeins}
\usepackage{multicol}
\usepackage{caption}
\usepackage{subcaption}

\usepackage{graphicx}
\graphicspath{{./}}
\usepackage[export]{adjustbox}

\usepackage{pgfplots}
\usepackage{pgfplotstable}
\usepackage{pgf-umlcd}
\usepackage{forest}
\usetikzlibrary{er,positioning,shapes,external}
\usepgfplotslibrary{colorbrewer,groupplots,external}

\pgfplotsset{compat=1.18}
\pgfplotsset{
    cycle list/Dark2-7,
    cycle multiindex* list={
        mark list*\nextlist
        Dark2-7\nextlist
    },
    every axis plot/.append style={thick}
}

\usepackage{xcolor}
\usepackage[normalem]{ulem}
\usepackage{comment}
\usepackage{todonotes}

\usepackage{pifont}

\usepackage{siunitx}
\theoremstyle{definition}

\newcommand{\uu}[1]{\mathbf{#1}}
\newcommand{\mean}[1]{\{\!\!\{#1\}\!\!\}}
\newcommand{\jump}[1]{[\![#1]\!]}

\newcommand{\ud}{\,\mathrm{d}}

\tikzset{
  bplus/.style={
    rectangle split,
    rectangle split horizontal,
    text width=1em,
    text centered,
    inner xsep=2pt,
    draw
  }
}

\forestset{
  multi edge/.style={
    edge path={
      \noexpand\path[\forestoption{edge}]
      (!u.#1 south)--(.child anchor)\forestoption{edge label};
    }
  }
}

\ExplSyntaxOn
\seq_new:N \l_tree_node_seq
\seq_set_from_clist:Nn \c_node_names_seq {one,two,three,four}
\NewDocumentCommand{\mpn}{m}{
  \seq_set_from_clist:Nn \l_tree_node_seq {#1}
  \seq_map_indexed_inline:Nn \l_tree_node_seq {
    \nodepart{\seq_item:Nn \c_node_names_seq {##1}} {##2}
  }
}
\ExplSyntaxOff

\usepackage[colorlinks=true]{hyperref}
\hypersetup{
    linkcolor=blue,
    filecolor=magenta,
    urlcolor=cyan
}

\begin{document}
\captionsetup[table]{skip=6pt}

\begin{frontmatter}
    \title{An agglomeration-based multigrid solver for the discontinuous Galerkin discretization of cardiac electrophysiology}

    \author[inst1]{Marco Feder}\corref{cor1}\ead{marco.feder@dm.unipi.it}
    \cortext[cor1]{Corresponding author.}

    \author[inst2]{Pasquale Claudio Africa}\ead{pafrica@sissa.it}

    \affiliation[inst1]{organization={Numerical Analysis Group, University of Pisa},
        addressline={Largo B. Pontecorvo, 5},
        city={Pisa},
        postcode={56126},
        country={Italy}}

    \affiliation[inst2]{organization={mathLab, International School for Advanced Studies (SISSA)},
        addressline={Via Bonomea, 265},
        city={Trieste},
        postcode={34136},
        country={Italy}}

    \begin{abstract}
        This work presents a novel agglomeration-based multilevel preconditioner designed to accelerate the convergence of iterative solvers for linear systems arising from the discontinuous Galerkin discretization of the monodomain model in cardiac electrophysiology. The proposed approach exploits general polytopic grids at coarser levels, obtained through the agglomeration of elements from an initial, potentially fine, mesh. By leveraging a robust and efficient agglomeration strategy, we construct a nested hierarchy of grids suitable for multilevel solver frameworks. The effectiveness and performance of the methodology are assessed through a series of numerical experiments on two- and three-dimensional domains, involving different ionic models and realistic unstructured geometries. The results demonstrate strong solver effectiveness and favorable scalability with respect to both the polynomial degree of the discretization and the number of levels selected in the multigrid preconditioner.
    \end{abstract}

\end{frontmatter}

\section{Introduction}\label{sec:intro}

Computational modeling of the heart has been actively pursued as a tool for accelerating cardiovascular research. However, the clinical applicability of such models is constrained by the high complexity and computational cost, for instance
when moving towards whole heart modeling and coupling different physics and scales. Hence, the development of numerical methods able to reduce the computing time while keeping accuracy becomes essential for
speeding-up fundamental research and, ultimately, for translation of modeling into clinical practice~\cite{Franzone2014MathematicalCE,Quarteroni_Manzoni_Vergara_2017}.

In this work, we focus specifically on cardiac electrophysiology models describing the electrical
activation of the myocardium. The electrical activity of the heart is
usually modeled through the bidomain or monodomain equations, which describe the
propagation of the transmembrane potential across the heart tissue. Both models consist
of an elliptic-parabolic system of nonlinear reaction-diffusion partial differential equations (PDEs). The
mathematical analysis of the system has already been addressed in~\cite{Franzone2014MathematicalCE}, and a wide strand of discretization approaches have been proposed
during the years, for instance in~\cite{COLLIFRANZONE200535,AFRICA2023111984,ANTONIETTI2025117957,ParallelNewtonKrylovBDDCFETIDP} and references therein. In order to take
into account the electrochemical reactions that occur at a cellular level, these systems are coupled, through
a non-linear reaction term, to a system of Ordinary Differential Equations (ODEs), modelling the inward and outward flow of ionic currents across the cell membrane~\cite{Franzone2014MathematicalCE}. During the years, a large variety of models have been developed, ranging from
reduced models with only one or few unknowns such as the Rogers-McCulloch ionic model~\cite{RogersMcCulloch}, or the FitzHugh and Nagumo model~\cite{FITZHUGH1961445}, up to the
Bueno-Orovio and ten Tusscher-Panfilov ionic models~\cite{TusscherPanfilov,TusscherPanfilov2,BUENOOROVIO2008544}. Due to the quick upstroke of the action potential, which is caused by voltage-dependent sodium channels, a numerically robust
calculation of the propagation of the wave across the tissue is well known to be computationally
challenging; indeed, the rapid increase of the transmembrane potential over a few milliseconds results in a steep wave front in space, requiring high resolution in both temporal and spatial
discretizations.

Overall, the resulting large-scale linear systems are prohibitive for direct solvers. In addition, the system matrix
often turns out to be very ill-conditioned, leading to a deterioration in the convergence rate of standard iterative methods. Hence, the main challenge to be addressed is
the reduction of the overall time-to-solution of the iterative method through the design of effective preconditioners~\cite{Benzi2005}.

The development of suitable preconditioners for such models is a vast and active area of research, and many approaches
have been proposed over the years. Among them, we mention domain-decomposition techniques such as Balancing Domain Decomposition by Constraints (BDDC)~\cite{BDDC1} and Finite
Element Tearing and Interconnecting (FETI) methods~\cite{FETI-DP}. Extensions of these methodologies, along with their
application to electrophysiology problems, have been extensively studied~\cite{ZampiniDualPrimalBidomain,ABDELHAMID2025118366,ParallelNewtonKrylovBDDCFETIDP,PavarinoScacchiAdditiveSchwarz}. A numerical study of AMG preconditioning
applied to the bidomain model can be found in~\cite{AMG_bidomain}.

Traditional strategies exploit linear continuous elements on very fine computational grids, in
order to capture the physical features of the solution. To cope with such complexity, higher order continuous elements or discontinuous Galerkin (DG) approaches have
recently received more attention~\cite{AFRICA2023111984,MonodomainHDG,botta2024highorderdiscontinuousgalerkinmethods}, thanks to their ability to capture
sharp gradients. For these reasons, in this work we will consider a DG discretization of the monodomain model. Our main focus is the development
of a preconditioner for accelerating the convergence of the large-scale linear system of equations.

For linear systems stemming from the discretization of elliptic operators, multigrid
methods are widely recognized as one of the most efficient approaches that scale to large parallel computers. Standard geometric
multigrid techniques require the construction of a hierarchy of meshes that is most often obtained by
(either uniform or adaptive) refinement of a coarse initial mesh. However, in the case of very complex geometries, such
as the ones stemming from realistic models, building such hierarchies is
far from trivial, requiring to fallback to non-nested methods~\cite{NonNested,MGNonNested,chen2024multigrid} or algebraic multigrid (AMG) solvers~\cite{Trottenberg2001}.

The use of polytopal shapes is attractive in this respect, as coarse grids can be simply generated by
agglomerating mesh elements together and various multilevel approaches have been proposed in the literature~\cite{Chan1998AnAM,BASSI201245,PAN2022110775,BOTTI2017382,MAVRIPLIS,DARGAVILLE2021113379}. However, providing automated
and flexible agglomeration strategies for polyhedral elements remains a challenging task which is subject of intense research~\cite{ANTONIETTI2026335,antonietti2025magnet,FEDER2025113773}.
Motivated by the good results observed in~\cite{FEDER2025113773}, we adopt the agglomeration algorithm developed therein to produce nested sequences of
agglomerated grids by exploiting R-trees~\cite{Guttman}, a spatial indexing data structure which excels at organizing geometric data
through bounding boxes, particularly in contexts where performing spatial queries for large sets of geometric objects in a fast way is required. This approach is fully
automated, robust, and dimension-independent. It is designed to be independent of the specific shape of the underlying
mesh elements, which can indeed be hexahedral, simplicial, polytopic, and stem from unstructured geometries. Notably, such
procedure is purely geometrical in that it depends exclusively on the initial mesh.

During the last few years, several multilevel strategies for methods posed on arbitrarily shaped elements have been investigated, such as in
Virtual Element methods (VEM)~\cite{AntoniettiVEMMG,AntoniettiVEMpMG,PradaPennacchioVEMAMG,antonietti2026reducedbasismultigridscheme} and Hybrid
High-Order (HHO) methods~\cite{di_pietro_high_order_mg_2023,HHOnonnestedMG}. For DG discretizations of Poisson problems posed on polytopal grids~\cite{cangiani2022hp}, both nested~\cite{hpMGAntonietti,antonietti2020agglomeration} and non-nested~\cite{AntoniettiPennesi} multigrid variants have
been proposed and theoretically analyzed. Polytopal DG methods have also been recently applied in the
field of brain electrophysiology~\cite{BLEIMERSAGLIO2025118249,LEIMERSAGLIO2026112}, which shares
the same abstract formulation of the model considered here~\cite{Schreiner2022}, and numerical investigations of non-overlapping Schwarz
preconditioners have been tested in~\cite{saglio2025massivelyparallelnonoverlappingschwarz}.

In this work, we propose a novel agglomeration-based multilevel preconditioner for the DG discretization of the monodomain model which exploits
the flexibility given by polytopic elements in creating hierarchies of coarse grids, starting from given unstructured meshes of interest. To increase
the computational efficiency of the proposed preconditioner, we combine our efficient agglomeration strategy with state-of-the-art matrix-free operator
evaluation kernels on the finest level of the hierarchy, thanks to the tensor-product structure of the basis functions~\cite{MatrixFreeDG}. We compare
the performance of our methodology with well-established preconditioners such as AMG and Block-Jacobi from the Trilinos library~\cite{Trilinos}, showing the effectiveness of our preconditioner in terms
of iteration counts of the preconditioned conjugate-gradient (PCG) solver and global time-to-solution. We validate our approach through several tests in
two and three dimensions, employing different ionic models, polynomial degrees, and geometries of different complexity. AMG and Block-Jacobi are well known
to be effective preconditioners and are a standard choice for discretizations with linear elements~\cite{FranzonePavarinoParallelMonodomain}. Our numerical
results indicate the advantage of our approach in terms of both iteration counts and time-to-solution, especially for higher polynomial degrees. In addition, we assess
the scalability of our preconditioner in a parallel computing environment, demonstrating its parallel efficiency.
The experiments are carried out using the \textsc{C++} software project \textsc{polyDEAL}~\cite{polydeal}, based on the \textsc{deal.II} finite element library~\cite{dealIIdesign}. It provides
building blocks for polytopal discontinuous DG methods and agglomeration-based multigrid, using Message Passing Interface (MPI) for distributed-memory parallelism.

The paper is structured as follows. In Section~\ref{sec:monodomain} we introduce the monodomain model for cardiac electrophysiology, and address its spatial and time discretization in Section~\ref{sec:discretization}. In Section~\ref{sec:rtree} we recall the R-tree data structure and the algorithmic realization of our approach, which is used to build hierarchies of agglomerated grids. We discuss the
details and setup of our preconditioner in Section~\ref{sec:prec_monodomain}. Finally, we present
several numerical results in Section~\ref{sec:numerical_experiments_monodomain}, and point to future developments in Section~\ref{sec:conc}.

\section[Mathematical model]{The mathematical model}\label{sec:monodomain}

\subsection{Mathematical modeling of cardiac electrophysiology}
The heart wall is organized into three layers: the thin inner endocardium, the thin outer epicardium, and the thick muscular middle
layer, the myocardium. The myocardium is mainly formed by cardiomyocytes, which are specialized, striated muscle cells that
are electrically excitable and drive the heart's mechanical function. When a cardiomyocyte receives an electrical stimulus, the electrochemical balance
across its membrane is altered, setting off a chain of biochemical events that change the cell's \emph{transmembrane potential} (the voltage difference
between the intracellular and extracellular environments). This produces a rapid
depolarization followed by a slower repolarization.
These voltage changes are governed by the opening and closing of voltage-gated ion channels, which selectively allow ions such as sodium,
potassium and calcium to cross the membrane.
Ionic currents change the membrane potential, while the membrane potential in turn controls
those currents. Neighboring cardiomyocytes are electrically coupled
by gap junctions, i.e. low-resistance intercellular channels, which permit the excitation to propagate from
cell to cell throughout the myocardium.

The mathematical description of these processes is made of two building blocks: a \emph{ionic model}, describing
the chemical processes taking place at the cellular scale, and an \emph{action potential propagation model}, representing the
spatial propagation of the action potential wavefront at the tissue scale. We will introduce both components in the following subsections.

\subsection{The monodomain model}
We start by fixing some notation. Given an open and bounded domain $D$, we denote by $L^2(D)$ the space of square integrable functions on $D$, endowed with the norm $\|\cdot\|_{0,D}$ induced by the inner
product $(\cdot,\cdot)_{D}$. We use the standard notation $H^s(D)$, with
associated norm $\|\cdot\|_{s,D}$, to indicate Sobolev spaces $W^{s,2}(D)$. Vector-valued functions are indicated in boldface.

In what follows we briefly introduce the monodomain model for cardiac electrophysiology. For a more comprehensive introduction to the topic, we refer the interested reader to~\cite{Franzone2014MathematicalCE,Quarteroni_Manzoni_Vergara_2017}.
Given an open, bounded domain $\Omega \subset \mathbb{R}^d$ ($d=2,3$) and a final time $T>0$, we introduce the transmembrane potential $u(\boldsymbol{x},t) \colon \Omega \times [0,T] \rightarrow \mathbb{R}$, and the
vector valued field $\boldsymbol{w}(\boldsymbol{x},t) \colon \Omega \times [0,T] \rightarrow \mathbb{R}^{S}$, where $S$ is the number of ionic variables involved in the ionic model at hand. The \emph{monodomain model} then
reads~\cite{Franzone2014MathematicalCE}:
\begin{equation}\label{eqn:monodomain}
    \begin{cases}
        \chi_m C_m \displaystyle\frac{\partial u}{\partial t}(t) - \nabla \cdot (\mathbb{D} \nabla u) + \chi_m \mathcal{I}_{\text{ion}}(u(t),\boldsymbol{w}(t)) = \mathcal{I}_{\text{app}}(\boldsymbol{x},t), & \text{in } \Omega \times (0,T],          \\[8pt]
        \displaystyle\frac{\partial \boldsymbol{w}(t)}{\partial t} = \boldsymbol{H}(u(t),\boldsymbol{w}(t)),                                                                                                  & \text{in } \Omega \times (0,T],          \\[8pt]
        \mathbb{D}\nabla u \cdot \boldsymbol{n} =0,                                                                                                                                                           & \text{on } \partial \Omega \times (0,T], \\[4pt]
        u(0)=u_0, \quad \boldsymbol{w}(0) = \boldsymbol{w}_0,                                                                                                                                                 & \text{in } \Omega.
    \end{cases}
\end{equation}
The conductivity tensor $\mathbb{D}$ is defined as: \[\mathbb{D} \coloneqq \sigma_l \boldsymbol{f}_0 \otimes \boldsymbol{f}_0 + \sigma_{t} \boldsymbol{s}_0 \otimes \boldsymbol{s}_0 + \sigma_{n} \boldsymbol{n}_0 \otimes \boldsymbol{n}_0, \] where
the vector fields $\boldsymbol{f}_0,\boldsymbol{s}_0,$ and $\boldsymbol{n}_0$ express the fiber, the sheetlet, and the sheet-normal (cross-fiber) directions, respectively. The longitudinal, transversal, and normal conductivities are denoted by $\sigma_l$, $\sigma_t$, and $\sigma_n$, respectively. Notice that $\boldsymbol{f}_0$, $\boldsymbol{s}_0, \boldsymbol{n}_0 \in [L^{\infty}(\Omega)]^3$, and $\sigma_l, \sigma_t, \sigma_n >0$ so that $\mathbb{D} \in [L^{\infty}(\Omega)]^{3 \times 3}$ is symmetric and positive definite almost anywhere in $\Omega$. The generation of myocardial
fibers has been performed using Laplace-Dirichlet Rule-Based Methods proposed and implemented in~\cite{PIERSANTI2021113468,Africa2023lifexfiber}.
The membrane cell capacitance is denoted by $C_m$, while $\chi_m$ is the membrane surface-to-volume ratio. The mathematical expression of the functions $\mathcal{I}_{\text{ion}}(u,\boldsymbol{w})$ and $\boldsymbol{H}(u,\boldsymbol{w})$ strictly depend on the choice of the ionic model, and will be specified later. The action potential is triggered by an external applied current $\mathcal{I}_{\text{app}}(\boldsymbol{x},t)$, which
mimics the presence of a (natural or artificial) pacemaker, while the $\mathcal{I}_{\text{ion}}(u,\boldsymbol{w})$ term
is responsible for describing the electric current generated by the flux of ionic species across the cell membrane.

We close the system with suitable initial conditions $u_0$, $\boldsymbol{w}_0$, and homogeneous Neumann boundary conditions on $\partial \Omega$ to model an electrically isolated domain. We denote by $\boldsymbol{n}$ the outward normal unit vector to the boundary.

\subsection{Ionic models}
Over the years, several ionic models have been proposed in the literature. Most of them are
written as systems of ordinary differential equations (ODEs):
\begin{equation}\label{eqn:generic_ode}
    \begin{cases}
        \frac{\partial \boldsymbol{w}(t)}{\partial t} = \boldsymbol{H}(u(t),\boldsymbol{w}(t)), & t \in (0,T], \\[8pt]
        u(0)=u_0, \quad \boldsymbol{w}(0) = \boldsymbol{w}_0,
    \end{cases}
\end{equation}where the unknowns are the transmembrane potential $u=u(t)$ and the vector $\boldsymbol{w}=(w_1,\ldots,w_S)$ collecting $S$ ionic variables. The dynamic
of the ionic variables is governed by the functions $\boldsymbol{H}(u,\boldsymbol{w})$ and $\mathcal{I}_{\text{ion}}(u,\boldsymbol{w})$ which couple the gating variables with
the evolution of the action potential. In this work, we will investigate two phenomenological ionic models: the FitzHugh-Nagumo~\cite{FITZHUGH1961445} and Bueno-Orovio~\cite{BUENOOROVIO2008544} models, which
we briefly describe hereafter.

\subsubsection{FitzHugh-Nagumo}
The FitzHugh-Nagumo model is described by the following scalar ODE:
\begin{equation}\label{eqn:FitzHugNagumo}
    \begin{cases}
        \frac{\partial w(t)}{\partial t} = \varepsilon(u-\Gamma w) & t \in (0,T], \\
        u(0)=u_0, w(0) = w_0.
    \end{cases}
\end{equation}The $\mathcal{I}_{\text{ion}}(u,w)$ term associated with this model has the expression:
\[\mathcal{I}_{\text{ion}}(u(t),w(t)) = k u(u-a)(u-1)+w,\]
where $(k,a,\varepsilon,\Gamma)$ are known parameters to tune the ionic model. Given its simplicity, it will
be used in the two-dimensional numerical experiments as a way to validate the proposed preconditioner.

\subsubsection{Bueno-Orovio}
The Bueno-Orovio model introduced in~\cite{BUENOOROVIO2008544} has $S\!=\!3$ gating variables $\boldsymbol{w}=(w_0,w_1,w_2)$ and is described by the following system of ODEs:
\begin{equation}\label{eqn:BuenoOrovio}
    \begin{cases}
        \frac{\partial w_0(t)}{\partial t} = [(b_0(u)-a_0(u))w_0 + a_0(u)w_0^{\infty}(u)] & t \in (0,T], \\
        \frac{\partial w_1(t)}{\partial t} = [(b_1(u)-a_1(u))w_1 + a_1(u)w_1^{\infty}(u)] & t \in (0,T], \\
        \frac{\partial w_2(t)}{\partial t} = [(b_2(u)-a_2(u))w_2 + a_2(u)w_2^{\infty}(u)] & t \in (0,T], \\
        u(0)=u_0, \boldsymbol{w}(0) = \boldsymbol{w}_0,
    \end{cases}
\end{equation}where
\begin{align*}
    a_0(u)          & =\frac{1 - H_{V_1}(u)}{H_{V_1}^{-}(u)(\tau_1^{''} - \tau_{1}^{'})+\tau_{1}^{'}},                            &
    a_1(u)          & =\frac{1 - H_{V_2}(u)}{H_{V_2}^{-}(u)(\tau_2^{''} - \tau_{2}^{'})+\tau_{2}^{'}},                              \\
    a_2(u)          & =\frac{1}{H_{V_2}(u)(\tau_3^{''} - \tau_3^{'})+\tau_3^{'}},                                                 &
    b_0(u)          & =-\frac{H_{V_1}(u)}{\tau_1^{+}},                                                                              \\
    b_1(u)          & =-\frac{H_{V_2}(u)}{\tau_2^{+}},                                                                            &
    b_2(u)          & =0,                                                                                                           \\
    w_0^{\infty}(u) & =1-H_{V_1^{-}}(u),                                                                                          &
    w_1^{\infty}(u) & =H_{V_o}(u) \Bigl(w_{\infty}^{\star} -1 + \frac{u}{\tau_2^{\infty}} \Bigr) + 1 - \frac{u}{\tau_2^{\infty}},   \\
    w_2^{\infty}(u) & =H_{V_3^{-}}^{K_3}(u).
\end{align*}
The function $H_{z_0}^{\varepsilon}(z) \coloneqq \frac{1+\tanh(\varepsilon(z-z_0))}{2}$ describes a smooth approximation\footnote{When $\varepsilon$ is omitted, it corresponds to the classical Heaviside function.} of the Heaviside function depending on the constant parameter $\varepsilon \in \mathbb{R}^{+}$. With this model, the ionic term is given by:
\begin{equation}
    \mathcal{I}_{\text{ion}}(u(t),\boldsymbol{w}(t))= \sum_{q\in \{fi,so,si\}} I^q(u(t),\boldsymbol{w}(t)),
\end{equation}where
\begin{align*}
    I^{fi} & = -\frac{H_{V_1}(u) (u-V_1)(\hat{V}-u)}{\tau_{fi}}w_0,                                                                                                                              &
    I^{si} & = -\frac{H_{V_2}(u)}{\tau_{si}} w_1 w_2,                                                                                                                                              \\
    I^{so} & = \frac{\left(1-H_{V_2}(u)\right) (u-V_o)}{H_{V_o}(u)(\tau_{o}^{''} - \tau_o^{'}) + \tau_o^{'}} + \frac{H_{V_2(u)}}{H_{V_{so}}(u)(\tau_{so}^{''} - \tau_{so}^{'}) + \tau_{so}^{'}}.
\end{align*}
Albeit relatively simple, the Bueno-Orovio model has the characteristic of capturing the main features of the electrophysiology in healthy myocardial tissues~\cite{GERBI2019,botta2024highorderdiscontinuousgalerkinmethods}. The parameter values used in
the numerical experiments are listed in~\ref{app:parameters}.

\section{Space and time discretization}\label{sec:discretization}
We can now introduce the weak formulation of Problem~\eqref{eqn:monodomain}. To this aim, we consider the Sobolev space $V\coloneqq H^1(\Omega)$, endowed with the classical $H^1$-norm $\|\cdot\|_{1,\Omega}$. We assume
that forcing terms, physical parameters, and initial conditions are sufficiently regular. The weak formulation then reads: find $u(t) \in V$, $\boldsymbol{w}(t) \in [V]^S$ such that for every $t \in (0,T]$:
\begin{equation*}\label{eqn:weak_monodomain}
    \begin{cases}
        \chi_m C_m \Bigl(\frac{\partial u(t)}{\partial t}, v \Bigr)_{\Omega} + a(u(t),v) + \chi_m (\mathcal{I}_{\text{ion}}(u(t),\boldsymbol{w}(t)),v)_{\Omega} = (\mathcal{I}_{\text{app}}(\boldsymbol{x},t),v)_{\Omega} & \forall v \in V,                     \\[8pt]
        \Bigl(\frac{\partial \boldsymbol{w}(t)}{\partial t},\boldsymbol{\psi}\Bigr)_{\Omega}  = (\boldsymbol{H}(u(t),\boldsymbol{w}(t)), \boldsymbol{\psi})_{\Omega}                                                      & \forall \boldsymbol{\psi} \in [V]^S, \\[8pt]
        u(0)=u_0, \quad \boldsymbol{w}(0) = \boldsymbol{w}_0                                                                                                                                                              & \text{in } \Omega,
    \end{cases}
\end{equation*}where the bilinear form $a(\cdot,\cdot)$ is defined as
\begin{equation}\label{eqn:bilinear_form}
    a(u,v) = \int_\Omega \mathbb{D} \nabla u \cdot \nabla v \ud \boldsymbol{x}.
\end{equation}Notice that the right-hand side of the equation involves only the external current term since homogeneous Neumann boundary conditions are imposed on $\partial \Omega$.

The spatial discretization of Problem~\eqref{eqn:monodomain} is performed with a standard
discontinuous Galerkin formulation~\cite{Riviere}. Let $\mathcal{T}_h$ be
a shape-regular triangulation of $\Omega$ consisting of disjoint elements $K$. We will assume that $\mathcal{T}_h$ is either a quadrilateral or hexahedral mesh. Interfaces are defined as the intersection of the
$(d-1)$-dimensional facets of neighboring elements. Further, we denote by $\mathcal{F}_h^I$ the union of all interior faces contained within $\Omega$ and by $\mathcal{F}_h^B$ those lying on the boundary $\partial \Omega$.
We define $\mathcal{Q}_p(K)$ as the space of tensor-product polynomials over the element $K$ of degree at most $p \geq 1$ in each coordinate direction, and the discontinuous finite element space over the mesh $\mathcal{T}_h$ as:

\[V_h^p = \bigl\{ v_h \in L^2(\Omega)\colon v_{h|K} \in \mathcal{Q}_p(K), K \in \mathcal{T}_h \bigr\}.\]Let $F \in \mathcal{F}_h^I$ be the common face between two neighboring elements $K^+$ and $K^-$, and
let $\boldsymbol{n}^+$ and $\boldsymbol{n}^-$ be the outward unit normal vectors, respectively. For regular enough scalar-valued function $v$ and a vector-valued function $\boldsymbol{q}$, we define
the average and jump operators as follows~\cite{UnifiedDG}:
\begin{align*}
    \mean{v}              & = \frac{1}{2}(v^+ + v^-),                           & \jump{v}              & = v^+\boldsymbol{n}^+ + v^-\boldsymbol{n}^-,                                         &  & \text{on } F \in \mathcal{F}_h^I, \\
    \mean{\boldsymbol{q}} & = \frac{1}{2}(\boldsymbol{q}^+ + \boldsymbol{q}^-), & \jump{\boldsymbol{q}} & = \boldsymbol{q}^+ \cdot \boldsymbol{n}^+ + \boldsymbol{q}^- \cdot \boldsymbol{n}^-, &  & \text{on } F \in \mathcal{F}_h^I.
\end{align*}The superscripts $\pm$ indicate the traces taken from the interior of elements $K^{\pm}$. This setting enables us to introduce the
following discrete bilinear forms $a_h(\cdot,\cdot): V_h^p \times V_h^p \rightarrow \mathbb{R}$ and $m_h(\cdot,\cdot): V_h^p \times V_h^p \rightarrow \mathbb{R}$:
\begin{align*}
    a_h(u_h,v_h) & = \int_\Omega \mathbb{D} \nabla_h u_h \cdot \nabla_h v_h \ud \uu{x}
    - \sum_{F \in \mathcal{F}_h^I}\int_{F} \Bigl( \mean{ \mathbb{D} \nabla u_h} \cdot \jump{v_h} + \mean{\mathbb{D} \nabla v_h} \cdot \jump{u_h}
    + \sigma \jump{u_h} \cdot \jump{v_h} \Bigr) \ud s,                                 \\
    m_h(u_h,v_h) & = \int_\Omega u_h v_h \ud \uu{x},
\end{align*}where $\nabla_h$ denotes the broken gradient operator defined element-wise.
The penalty parameter is defined edge-wise as $\sigma(\boldsymbol{x})\!=\!\alpha\boldsymbol{n}^T \mathbb{D} \boldsymbol{n} (\boldsymbol{x})>0$, with $\alpha=\frac{p(p+1)|F|_{d-1}}{|K|_d}$, being $|\cdot|_{d}$ the $d$-dimensional Hausdorff measure of the corresponding geometric entity.

\subsection{Semidiscrete and fully-discrete formulation}
The semi-discrete formulation of Problem~\eqref{eqn:monodomain} then reads: for any $t \in(0,T]$, find $u_h(t) \in V_h^p$ and $\boldsymbol{w}_h(t) \in [V_h^p]^S$ such that
\begin{equation*}\label{eqn:semidiscrete}
    \begin{cases}
        \chi_m C_m\Bigl(\displaystyle\frac{\partial u_h(t)}{\partial t},v_h \Bigr)_{\Omega} + a_h(u_h(t),v_h) + \chi_m\Bigl(\mathcal{I}_{\text{ion}}(u_h(t),\boldsymbol{w}_h(t)),v_h\Bigr)_{\Omega} = \Bigl(\mathcal{I}_{\text{app}}(\boldsymbol{x},t),v_h\Bigr)_{\Omega} & \forall v_h \in V_h^p, \\[8pt]
        \displaystyle\frac{\partial \boldsymbol{w}_h(t)}{\partial t} = \boldsymbol{H}(u_h(t),\boldsymbol{w}_h(t))                                                                                                                                                      ,                           \\[8pt]
        u_h(0)=u_{h,0}, \quad \boldsymbol{w}_h(0)=\boldsymbol{w}_{h,0}                                                                                                                                                                                                    & \text{in } \Omega,
    \end{cases}
\end{equation*}where $u_{h,0}$ and $\boldsymbol{w}_{h,0}$ are given initial conditions. Let $N$ be the dimension of the discrete space $V_h^p$, and let $\{\phi_i\}_{i=1}^N$ be
the basis functions of $V_h^p$ based on tensor products of Lagrangian polynomials. Then, for all $l=1,\ldots,S$, we have the expansions $u_h(t)\!=\sum_{i=1}^N U_i(t) \phi_i$ and $w_{h}^{l}(t)\!=\sum_{i=1}^N W_{h,i}^l(t) \phi_i$. Here $\boldsymbol{U}(t) \in \mathbb{R}^N$
and $\boldsymbol{W}_{h,l}(t) \in \mathbb{R}^N$ are vectors holding the expansion coefficients of $u_h(t)$ and $w_{h}^{l}(t)$ with respect
to the basis $\{\phi_i\}_{i=1}^N$ for each time instant $t\in (0,T]$.
We can hence define the relevant matrices, forcing terms and vectors as follows:
\begin{eqnarray}
    \begin{aligned}\label{eqn:matrices}
        (A)_{i,j}                                          & = a_h(\phi_j,\phi_i),                                           & \qquad
        (M)_{i,j}                                          & = m_h(\phi_j,\phi_i),                                                    \\
        \bigl(I_{\text{ion}}^h(u,\boldsymbol{w})\bigr)_{i} & = (\mathcal{I}_{\text{ion}}(u,\boldsymbol{w}),\phi_i)_{\Omega}, & \qquad
        \bigl(I_{\text{app}}^h\bigr)_{i}                   & = (\mathcal{I}_{\text{app}},\phi_i)_{\Omega}.                   & \qquad
    \end{aligned}
\end{eqnarray}

With this notation, the algebraic formulation of the semi-discrete problem reads:
\begin{equation}\label{eqn:semidiscrete_algebraic}
    \begin{cases}
        \chi_m C_m M \dot{\boldsymbol{U}}_h + A\boldsymbol{U}_h  + \chi_m I_{\text{ion}}^h(\boldsymbol{U}_h,\boldsymbol{W}_h) = I_{\text{app}}^h(t), \\
        \dot{\boldsymbol{W}}_h = \boldsymbol{H}(\boldsymbol{U}_h,\boldsymbol{W}_h),                                                                  \\
        \boldsymbol{U}_h(0) = \boldsymbol{U}_{0},                                                                                                    \\
        \boldsymbol{W}_h(0) = \boldsymbol{W}_{0}.
    \end{cases}
\end{equation}To obtain the fully-discrete system, we partition the time interval $[0,T]$ into $N_T$ uniform subintervals, each of length $\Delta t\!=\!t_{n+1}-t_n$, with $t_n=n \Delta t$ for $n=0,\ldots,N_{T}-1$. We use the subscript $n$ to indicate
the approximation of $\boldsymbol{U}_h(t)$ and $\boldsymbol{W}_h(t)$ at time $t_n$. For the temporal discretization, we adopt a second-order Backward Differentiation Formula (BDF2). Given initial conditions $\boldsymbol{U}_0$ and $\boldsymbol{W}_0$, together with
suitable initializations of $\boldsymbol{U}_{1}$ and $\boldsymbol{W}_{1}$, the discrete scheme is: for $n=1,\ldots,N_T-1$, find $\boldsymbol{U}_{n+1}$ and $\boldsymbol{W}_{n+1}$ such that
\begin{equation}
    \chi_m C_m M\frac{3\boldsymbol{U}_{n+1}-4\boldsymbol{U}_n+\boldsymbol{U}_{n-1}}{2\Delta t} + A\boldsymbol{U}_{n+1} + \chi_m I_{\text{ion}}^h(\boldsymbol{U}^{*},\boldsymbol{W}_{n+1}) =  I_{\text{app},n+1}^h,
\end{equation}which yields an explicit formula for $\boldsymbol{U}_{n+1}$
\begin{equation}\label{eqn:fully_discrete}
    \Bigg(3\chi_m C_m\frac{M}{2\Delta t} + A \Bigg)\boldsymbol{U}_{n+1}= I_{\text{app},n+1}^h -  \chi_m I_{\text{ion}}^h(\boldsymbol{U}^{*},\boldsymbol{W}_{n+1}) + \chi_m C_m \frac{M}{2\Delta t} (4\boldsymbol{U}_n - \boldsymbol{U}_{n-1}).
\end{equation}Setting $\boldsymbol{U}^* = \boldsymbol{U}_{n+1}$ yields a fully implicit BDF2 scheme. To overcome this issue, we use the extrapolation formula $\boldsymbol{U}^* = 2\boldsymbol{U}_n - \boldsymbol{U}_{n-1}$, which
is second-order accurate with respect to $\Delta t$~\cite{GERVASIO2006347,AFRICA2023111984}. Therefore, given the solution $\bigl(\boldsymbol{U}_n,\boldsymbol{W}_n\bigr)$ at time $t_n$, the solution at time
$t_{n+1}$ is computed as follows:
\begin{itemize}
    \item Solve the ionic model for $\boldsymbol{W}_{n+1}$ at each degree of freedom,

    \item Plug the computed $\boldsymbol{W}_{n+1}$ in $I_{\text{ion}}(\boldsymbol{U}^{*},\boldsymbol{W}_{n+1})$ into Equation~\eqref{eqn:fully_discrete},

    \item Solve the resulting linear system for $\boldsymbol{U}_{n+1}$.
\end{itemize}
Note that this scheme allows crucial gains in computational efficiency: the problem for $\boldsymbol{U}_{n+1}$ is
linear and hence the system matrix $3\chi_m C_m \frac{M}{2\Delta t} + A$ and the relative preconditioner
are assembled once and for all at the beginning of the computation.

As for the evaluation of the ionic term $I_{\text{ion}}^h(\boldsymbol{U}^{*},\boldsymbol{W}_{n+1})$ in~\eqref{eqn:fully_discrete}, we opt
for the so-called Ionic Current Interpolation (ICI) method~\cite{ICI}: we evaluate the ionic term at the degrees of freedom, and then interpolate
onto quadrature nodes on mesh elements while assembling local contributions for the right-hand side vector. This procedure is cheaper and requires less memory than solving the system of ODEs at each quadrature node to compute $I_{\text{ion}}^h$ (a procedure also known as State Variable Interpolation). The system of ODEs for $\boldsymbol{W}_{n+1}$ is solved
using a BDF2 scheme.

Since the system matrix is symmetric and positive definite, the resulting linear system is solved with
the preconditioned conjugate gradient (PCG) method. The focus of the remainder of this work will be the development of a
multilevel preconditioner which exploits an efficient coarsening procedure of the given mesh to
accelerate the convergence of the PCG solver.

\section{R-tree based agglomeration}
\label{sec:rtree}
Before introducing our multigrid preconditioner, we briefly recall the basic properties of the R-tree data structure originally proposed
in~\cite{Guttman} and discuss its application to construct agglomerated grids following our previous work~\cite{FEDER2025113773}.

R-trees are hierarchical data structures used for the dynamic\footnote{Dynamic here means that no global reorganization is required upon insertion or deletion of new elements of the tree.} organization of
collections of $d$-dimensional geometric objects, representing them by their \emph{Minimum Bounding -- axis aligned -- d-dimensional Rectangle}, also denoted by MBR. The minimum bounding rectangle is often
referred to as \emph{bounding box}. An internal node of an R-tree consists of the MBR that bounds all its children. In particular, each internal node stores two data: a way of identifying a
child node and the MBR of all entries within this child node. The actual geometric data are stored in
the \emph{leafs} of the tree, that is, in its terminal nodes.

The original R-tree is based solely on the minimization of the measure of each MBR. Several variants
have been proposed, aimed at either improving performance or flexibility, generally depending on the domain
of application. It is often desirable to minimize the \emph{overlap} between MBRs. Indeed, the larger the overlap, the larger the number of paths to
be processed during queries. Moreover, the smaller the overlap the closer the agglomerated element will conform
to the corresponding bounding boxes, thus making the resulting agglomerated grid qualitatively
close to rectangular. Among the several available variants, we adopt the R$^*$-tree data structure designed in \cite{Beckmann1990TheRA}. In the remainder
of the paper, we will not make any distinction between R-trees and R$^*$-trees, as we always employ the latter.

In our implementation, we rely on the \textsc{Boost.Geometry} module supplied by the \textsc{Boost} C++ libraries~\cite{BoostLibrary} for the construction and manipulation of R-trees. Its kernels
are designed to be agnostic with respect to the number of space dimensions, coordinate systems, and
types.

\subsection{Agglomeration algorithm}
\label{subsec:rtree-agglomeration}
Here we describe the construction of an agglomerated mesh starting from a given finite element mesh, together with the associated R-tree data structure. Let $\mathcal{T}_h$ be a given mesh and
assume the associated R-tree of order $(m,M)$ has been constructed. We denote by $L$ the total number of levels of the R-tree and with
$N_l$ the set of nodes in level $l$, for $l \in \{0,\ldots,L-1\}$. Our agglomeration strategy depends on an input parameter $l \in \{0,\ldots,L-1\}$ which is
the index of the level to be used to generate the final agglomerates. The basic idea consists of looping through the nodes $N_l$ and, for each node $n \in N_l$, descending recursively its children
until leaf nodes are reached\footnote{Such technique of traversing a tree data structure is also known as depth-first search (DFS).}. These nodes share the same ancestor node $n$ (on level $l$) and thus are \emph{agglomerated}. We store such elements in a vector $v[n]$, indexed by the
node $n$. This procedure is outlined in Algorithm \ref{alg:r_tree_algo} and Algorithm~\ref{alg:recursive_extract}.

\begin{algorithm}[H]
    \caption{Creation of agglomerates.}
    \begin{algorithmic}[1]
        \Require{R-tree $R$ of order $(m,M)$}
        \Require{$l \in \{0, \ldots,L-1\}$ target level.}
        \Ensure{$v$ vector s.t. $v[n]$ stores leafs associated to node $n$}
        \Statex
        \Function{ComputeAgglomerates}{$R,l$}
        \For{node $n$ in $N_l$}
        \State $v[n] \gets $ \Call{Extractleaves}{$l,n$} \Comment{Defined in Algorithm~\ref{alg:recursive_extract}}
        \EndFor
        \State \Return $v$
        \EndFunction
    \end{algorithmic}
    \label{alg:r_tree_algo}
\end{algorithm}

\begin{algorithm}[H]
    \caption{Recursive extraction of leafs from a node $n$ on level $l$.}
    \begin{algorithmic}[1]
        \Require{$l \in \{0, \ldots,L-1\}$ target level}
        \Require{$n \in N_l$}
        \Ensure{vector $v$ containing leafs which share the ancestor node $n$.}
        \Statex
        \Function{ExtractLeafs}{$l,n$}
        \If{$l=0$}
        \State $v.\Call{push\_back}{\Call{get\_leafs}{n}}$
        \State \Return $v$
        \EndIf
        \For{node $ch_n$ in \Call{children}{$n$}}
        \State  $v.\Call{push\_back}{\Call{ExtractLeafs}{l-1,ch_n}}$
        \EndFor
        \State \Return $v$
        \EndFunction
    \end{algorithmic}
    \label{alg:recursive_extract}
\end{algorithm}
After the recursive visit of the children of a node $n$, a list of original mesh elements
can be stored in the vector $v[n]$ and flagged appropriately for agglomeration. Since the R-tree data structure provides a spatial partition of mesh elements, each of which is
uniquely associated with a node, the traversal of all nodes on a level $l$ provides a partition of mesh elements into agglomerates. The overall agglomeration procedure can be summarized as follows:
\begin{itemize}
    \item Compute $\{\text{MBR}(T_i)\}_i$ for $i=1,\ldots,\text{Card}(\mathcal{T}_h)$,
    \item Build the R-tree data structure using $\{\text{MBR}(T_i)\}_i$,
    \item Choose one level $l \in \{0, \ldots,L-1 \}$ and apply Algorithm \ref{alg:r_tree_algo},
    \item For each node $n$, flag together mesh elements stored in $v[n]$.
\end{itemize}

In practice, elements in $v[n]$ are mesh iterators, i.e. lightweight objects such as pointers that
uniquely identify elements of $\mathcal{T}_h$. The construction of the R-tree and the actual identification
of agglomerates turns out to be a highly efficient operation also for large three-dimensional geometries. The interested reader
is referred to~\cite{FEDER2025113773} for a detailed analysis of the computational cost of the proposed agglomeration strategy and its
application to several three-dimensional geometries. Notably, Algorithm~\ref{alg:r_tree_algo} can be employed to generate sequences of \emph{nested} agglomerated meshes for multilevel methods. As a consequence, such meshes can be used as levels of a hierarchy in a multigrid algorithm, allowing the usage
of simpler and much cheaper intergrid transfer operators, when compared to non-nested variants; as a matter of fact, efficient intergrid transfer operators for non-matching meshes
are far from trivial even for simple geometries and their construction constitutes a severe bottleneck in terms of computational efforts, becoming critical in 3D~\cite{AntoniettiPennesi,HHOnonnestedMG}. In the context of
Lagrangian finite elements on standard hexahedral or quadrilateral grids, a completely parallel and matrix-free implementation of the non-nested geometric multigrid method has been recently addressed in~\cite{MGNonNested}.

\begin{figure}[htbp]
    \centering

    \begin{subfigure}{0.6\textwidth}
        \centering
        \includegraphics[width=0.6\textwidth]{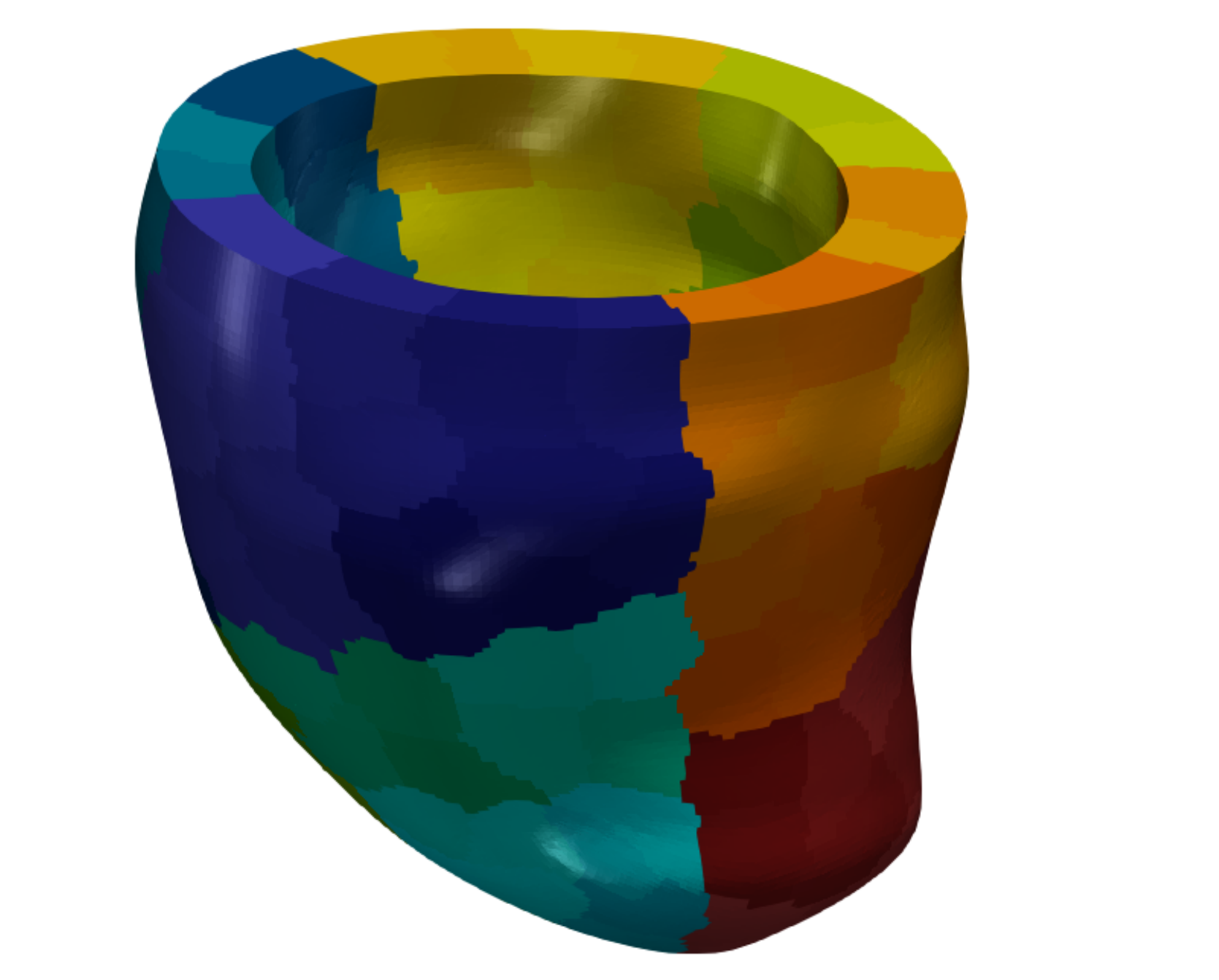}
        \caption{Hexahedral ventricle mesh partitioned among $128$ MPI processes.}
        \label{fig:ventricle_partition}
    \end{subfigure}

    \vspace{0.5cm}

    \begin{subfigure}{0.3\textwidth}
        \centering
        \includegraphics[width=0.8\textwidth]{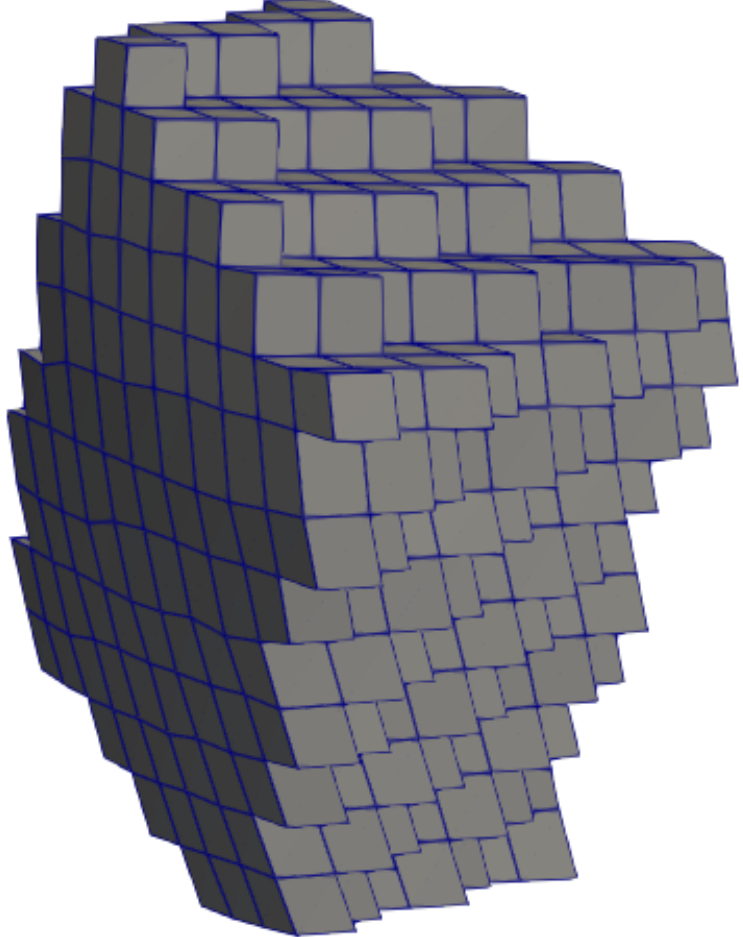}
        \caption{}
        \label{fig:ventricle_agglomerate}
    \end{subfigure}
    \hfill
    \begin{subfigure}{0.3\textwidth}
        \centering
        \includegraphics[width=0.8\textwidth]{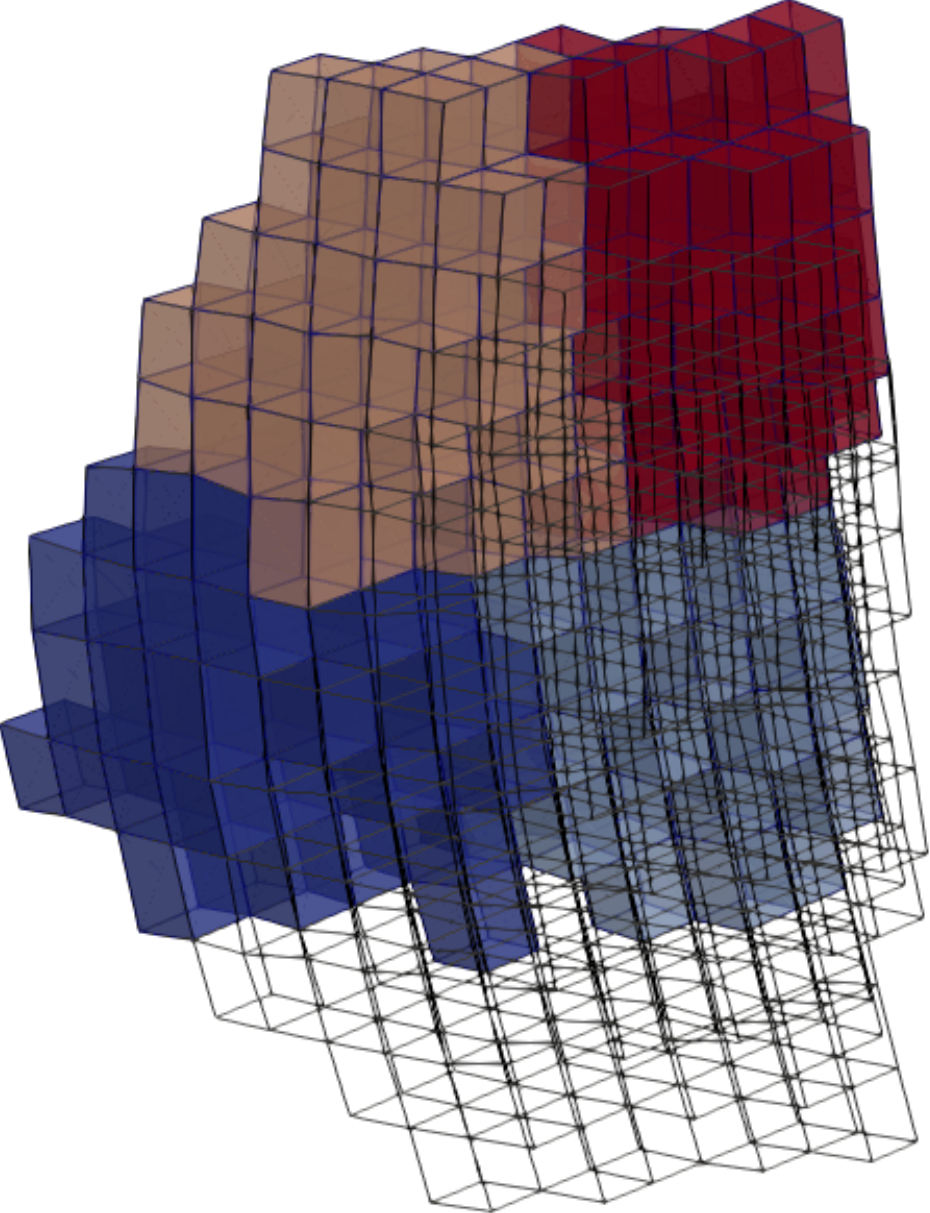}
        \caption{}
        \label{fig:ventricle_sub_agglo1}
    \end{subfigure}
    \hfill
    \begin{subfigure}{0.3\textwidth}
        \centering
        \includegraphics[width=0.8\textwidth]{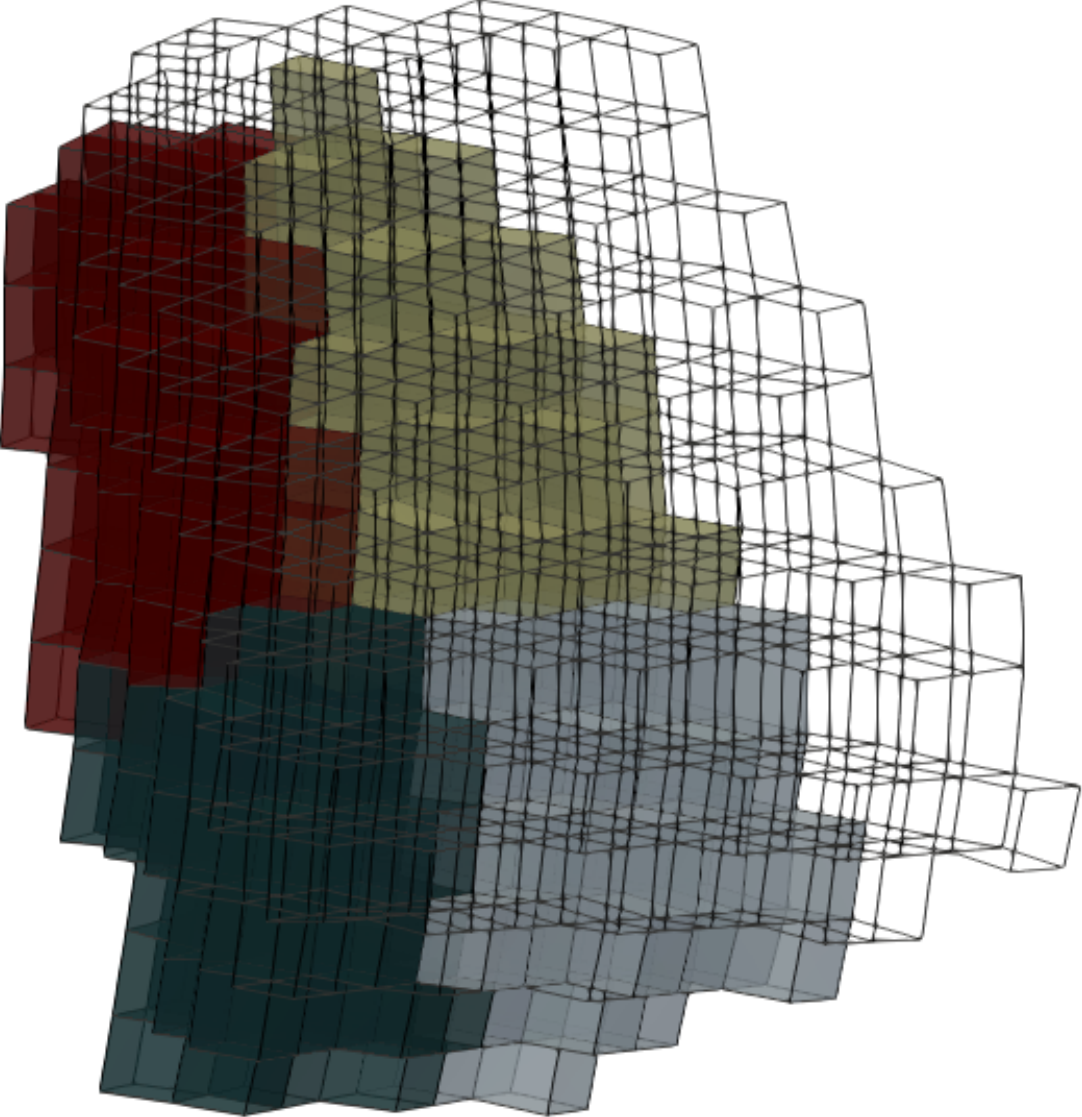}
        \caption{}
        \label{fig:ventricle_sub_agglo2}
    \end{subfigure}

    \caption{
        (\subref{fig:ventricle_partition}) MPI partitioning of the realistic hexahedral ventricle mesh into $128$ subdomains. Each color corresponds to a different MPI process.
        Bottom row: detailed views of an agglomerate on the boundary of the domain
        (\subref{fig:ventricle_agglomerate}) and its $8$ sub-agglomerates
        (\subref{fig:ventricle_sub_agglo1}) and
        (\subref{fig:ventricle_sub_agglo2}), each displayed in a different color.
    }
    \label{fig:ventricle_combined}
\end{figure}

We end this section by showing an example of agglomeration based on R-trees for a realistic 3D mesh of
a left ventricle, made of $374,022$ hexahedral elements, obtained from real personalized geometries~\cite{left_ventricle_mesh}. In practical simulations, the computational workload is
split among processes through a mesh partitioner such as \textsc{parmetis}~\cite{METIS}. Then, each MPI process owns a subdomain over which it performs local computations. A visualization of the MPI partitioning of the ventricle mesh into
$128$ subdomains is shown in Fig.~\ref{fig:ventricle_partition}. Next, \emph{within} each locally owned partition, each process generates agglomerates
using the tree data structure outlined above. An example agglomerate element (on the boundary of the ventricle mesh) is shown in Fig.~\ref{fig:ventricle_agglomerate}, in grey. Its $8$ sub-agglomerates, belonging
to the next finer level of the hierarchy, are depicted in Figs.~\ref{fig:ventricle_sub_agglo1} and~\ref{fig:ventricle_sub_agglo2}. For better
visualization, the sub-agglomerates, identified by different colors, are displayed in either of the two plots.
We do not show results with agglomeration based on graph partitioning, but we refer the reader to~\cite{FEDER2025113773} for several comparisons between
the two approaches. A detailed analysis of different agglomeration strategies in the context of multilevel solvers can be found in~\cite{DARGAVILLE2021113379}.

\section{Agglomeration-based multilevel preconditioning}\label{sec:prec_monodomain}
In this section, we propose an agglomeration-based multilevel preconditioner to accelerate the convergence of the PCG solver. Geometric multigrid methods
typically start from a coarse mesh which is then repeatedly subdivided in finer meshes to generate a hierarchy of levels for which intergrid transfers are simple and fast. Hence,
they hinder starting from a \emph{predetermined} fine grid.

A popular and effective way to obviate this issue is to use AMG methods, which require only
explicit knowledge of the system matrix, without relying on any geometrical information.
Algebraic multigrid methodologies are widely employed for linear systems arising from scalar elliptic
PDEs, and are renowned for their scalable performance~\cite{hypre,DAMBRAAMG}.

When applied to finite element discretizations, AMG preconditioners typically target matrices arising from linear
continuous elements. In the discontinuous Galerkin setting, however, the
redundancy of degrees of freedom associated with the same grid point is
known to hinder the creation of aggregates. To address this issue, AMG
solvers tailored to DG schemes were developed in~\cite{Bastian2012,AntoniettiAMGDG}.

Here, instead, we assume a target fine mesh $\mathcal{T}_h$ over which we need to solve the problem to be given. Rather
than uniformly refining $\mathcal{T}_h$, or feeding the system matrix to an AMG solver, we take
a hybrid route and use the built-in flexibility of DG methods, in conjunction with
the coarsening strategy developed in~\cite{FEDER2025113773}, to build a
sequence of \emph{nested} meshes and \emph{coarser} operators that can be employed as
level matrices in a multilevel method. The nestedness of the sublevels implies that
transfer operators are much cheaper to compute compared to non-nested variants.
All in all, the resulting multilevel preconditioner will be built in an algebraic manner, injecting geometric
information delivered by the agglomeration routine. Therefore, we propose a
method that is somewhat dual in that it can be interpreted both as algebraic and geometric.

The convergence properties of our approach rely on the theoretical framework developed in~\cite{AntoniettiSartiVeraniDGhpMG} for multigrid
methods applied to Interior Penalty DG discretizations on meshes of agglomerated elements.

\subsection{Inherited bilinear forms}
Using the agglomeration procedure described in Section~\ref{sec:rtree}, we build a sequence of $L$ nested meshes $\{\mathcal{T}_l\}_{l=0}^{L-1}$ with associated finite element
spaces $\{V_{l}\}_{l=0}^{L-1}$. By construction, the original mesh and finite element space correspond to $\mathcal{T}_0 \equiv \mathcal{T}_h$ and $V_0 \equiv V_h^p$, respectively. We index the sequence
of meshes and spaces so that to the lowest index corresponds the finest level. Since grids are nested, also the corresponding
finite dimensional subspaces are nested: $$V_{L-1} \subset V_{L-2} \subset \cdots \subset V_{0}.$$

On the finest grid $\mathcal{T}_0$, we consider the original bilinear form $\mathcal{A}_0: V_0 \times V_0 \rightarrow \mathbb{R}$ defined by:

\begin{equation}
    \mathcal{A}_0(u,v)\coloneqq 3\chi_m C_m \frac{m_h(u,v)}{2 \Delta t} + a_h(u,v)  \qquad \forall u,v \in V_{0},
\end{equation}with matrix representation $A_0\coloneqq 3\chi_m C_m \frac{M}{2\Delta t} + A$, as introduced in Section~\ref{sec:discretization}.

On each bounding box $B_K$ of an agglomerated polytopic element $K \in \mathcal{T}_l$, $l \in \{1, \ldots, L-1\}$, the basis functions are
constructed by defining on $B_K$ the standard polynomial space $\mathcal{Q}^p(B_K)$. As $\overline{K} \subset \overline{B_K}$, the basis on $K$ is defined
by restricting each basis function to $K$. Other choices are possible, and we refer to the textbook~\cite{BookDGPolytopes} for an extensive discussion and relevant implementation
details.

There are essentially two ways to build a sequence of coarser operators: the \emph{inherited} approach~\cite{AntoniettiSartiVeraniDGhpMG}, where the discrete operators are recursively built
by restricting the fine grid operators, and the \emph{non-inherited} one, where bilinear forms
are explicitly assembled on each agglomerated level $\{\mathcal{T}_{l}\}_{l=1}^{L-1}$ of the hierarchy.
Assembling bilinear forms on the coarser levels involves the computational burden of numerically integrating
over agglomerated elements.
To avoid such cost, we consider the inherited approach, where operators are defined iteratively from the
restriction of the original bilinear form $\mathcal{A}_0$ as follows:
\begin{equation}\label{eqn:inherited}
    \mathcal{A}_{l+1}(u,v) \coloneqq \mathcal{A}_{l}\bigl(\mathcal{P}_{l+1}^{l}u,\mathcal{P}_{l+1}^{l}v\bigr) \qquad \forall u,v \in V_{l+1},
\end{equation}with associated matrices:
\begin{equation}\label{eqn:inherited_matrices}
    A_{l+1} \coloneqq \bigl(\mathcal{P}_{l+1}^l\bigr)^T A_l \mathcal{P}_{l+1}^l,
\end{equation}where $\mathcal{P}_{l+1}^{l}$ denotes the prolongation from the coarser level
$l+1$ to the finer level $l$. As usual, the restriction operator $\mathcal{R}_{l}^{l+1}$ is defined as the adjoint of the prolongation. This technique for constructing coarser operators is one of the fundamental components of AMG methods and is commonly referred to as the \emph{Galerkin projection} or \emph{Galerkin triple product}. The inherited procedure is summarized in Algorithm~\ref{alg:agglomeration_inherited}.

From the implementation standpoint, we observe that each prolongation matrix $\mathcal{P}_{l+1}^l \in \mathbb{R}^{n_l \times n_{l+1}}$, where $n_l=\dim V_l$ and  $n_{l+1}=\dim V_{l+1}$, is a distributed matrix whose parallel layout is easily defined thanks to nestedness of the partitions and the
fact that in a DG discretization degrees of freedom are non-overlapping. Moreover, the simple block structure of the prolongation matrices allows to perform restriction and prolongation
of vectors in a matrix-free fashion by using the local action of each block element-wise.

Avoiding assembling the coarser operators $\{A_l\}_{l=1}^{L-1}$ explicitly is quite attractive, as looping over all the polytopes on every level becomes
increasingly expensive when agglomerated elements are composed of many sub-elements, particularly for high polynomial degrees, due to the large number of quadrature points. In principle, if homogeneous coefficients are considered, such cost could be dramatically reduced by using quadrature-free approaches for polytopic elements that reduce the integration over a polytope to only boundary evaluations~\cite{AntoniettiHouston}. Regardless of the quadrature
approach, the original mesh $\mathcal{T}_h$ needs to be exploited to assemble all operators, with a huge gain in favor of the quadrature-free approach if it can be employed.

Another appealing feature of the inherited technique is its easiness to be added in existing finite
element frameworks. In particular, its parallel implementation is simpler compared to the non-inherited approach. Indeed, in the non-inherited
approach, flux computations on processor partition boundaries require accessing data from ghost agglomerated elements, which are composed
of many layers of fine ghost cells that are not locally owned. With the inherited
approach, instead, the polytopic shapes obtained by agglomeration serve as a tool to build the sequence of
prolongation operators $\{\mathcal{P}_{l+1}^{l}\}_{l \geq 0}$, which are then used to
construct coarse-level matrices.

On the other hand, for multigrid schemes applied to DG discretizations of the Poisson problem on standard grids, it has been
proven that only the non-inherited multigrid provides uniform convergence with respect to the number of levels~\cite{AntoniettiSartiVeraniDGhpMG}.
Such behavior was verified also in~\cite{BOTTI2017382} when using one V-cycle of agglomerated multigrid as a preconditioner, and
subsequently fixed by adjusting the amount of stabilization in each level.

Remaining in the realm of agglomerated multigrid, we mention that in~\cite{FEDER2025113773} the
non-inherited variant applied to hierarchies generated by the agglomeration
algorithm in Section~\ref{sec:rtree} has been shown to work effectively, yielding
uniform convergence with respect to the number of levels for several test cases.
However, for the model considered in this work, numerical evidence shows that the number
of iterations is not affected by the number of levels also when the inherited approach is used, cf. Section~\ref{sec:numerical_experiments_monodomain} for details.

We point out that, despite this work focuses on the monodomain model, we foresee the application of the proposed preconditioner to a variety of other symmetric and positive definite problems.

\begin{algorithm}
    \caption{Construction of level operators $\{A_l\}_l$.}
    \label{alg:agglomeration_inherited}
    \begin{algorithmic}[1]
        \Require Mesh $\mathcal{T}_h$.
        \Ensure Sequence of level operators $\{A_l\}_l$.

        \Function{GenerateCoarserOperators}{}
        \State \Call{read\_mesh}{$\mathcal{T}_h$}
        \State $A_0 \gets \Call{assemble\_system}{\mathcal{T}_h}$
        \State \Call{build\_agglomerated\_hierarchy}{$\mathcal{T}_h$} \Comment{Coarser levels $\{\mathcal{T}_l\}_{l=1}^{L-1}$ (Section~\ref{sec:rtree})}
        \For{$l \gets 0$ to $L-1$}
        \State $\mathcal{P}_{l+1}^l \gets \Call{compute\_prolongation}{l,l+1}$
        \State $A_{l+1} \gets \bigl(\mathcal{P}_{l+1}^l\bigr)^T A_l \mathcal{P}_{l+1}^l$ \Comment{Galerkin product}
        \EndFor
        \State \Return $\{A_l\}_l$
        \EndFunction
    \end{algorithmic}
\end{algorithm}

\subsection{Matrix-free operator evaluation}
\label{subsec:matrix-free}
Within classical Krylov methods, the main computational kernel is the evaluation of the action of the operator on a vector. On current hardware architectures, this is
rooted in high memory traffic from loading matrix entries into compute units performing arithmetic operations. For polynomial degrees $p \geq 2$, the increased
bandwidth of the system matrix leads to a significant reduction in the throughput per degree of freedom. Matrix-free algorithms, conversely, avoid storing the matrix representation
of the operator in order to reduce memory traffic, typically yielding improved performance for $p \geq 2$.
Hence, whenever the iterative solver needs to evaluate
the action of the operator, the matrix-free approach evaluates the spatial integrals of the finite element discretization by numerical quadrature rules.
State-of-the-art matrix-free implementations exploit sum-factorization techniques for tensor-product basis functions~\cite{ORSZAG198070} and parallelization via mesh partitioning~\cite{p4est}. When tensor-product shape functions are combined with tensor-product quadrature rules, this computation becomes highly efficient through sum-factorization, which transforms the
standard matrix-based operator evaluation involving $(p+1)^{2d}$ operations for $(p+1)^d$ degrees of freedom with polynomial degree $p$ in $d$ spatial dimensions at $(p+1)^d$ quadrature points into a sequence of $d$ one-dimensional operations, each requiring $(p+1)^{d+1}$ operations. For three-dimensional problems, this approach reduces the algorithmic complexity by two orders in $p$. Additionally, SIMD
vectorization (Single Instruction Multiple Data) over batches of elements is employed to further enhance performance~\cite{KRONBICHLER2012135,MatrixFreeDG}.

Irrespective of the chosen preconditioner, the action of the finest operator $3\chi_m C_m \frac{M}{2\Delta t} + A$ within the PCG solver is computed using the matrix-free kernels available, while coarse-level operators $\{A_l\}_{l=1}^{L-1}$ are applied in a matrix-based way due to the lost of the
tensor product structure of the basis functions on agglomerated elements. With reference to our agglomeration-based preconditioner, also
the Chebyshev smoother exploits the efficient operator evaluation on the finest level. This hybrid approach exploits the reduced memory footprint of the matrix-free evaluation on the fine level, where the number of degrees of freedom is the largest, while
keeping the computational cost of coarse-level operator applications low, as the size of the coarse problems is increasingly smaller. In summary, within
the PCG solver, the action of the finest-level matrix $3\chi_m C_m \frac{M}{2\Delta t} + A$ on a vector is computed in a matrix-free fashion, while the
application of the preconditioners relies on assembled matrices on coarser levels.  The evaluation
of the right-hand side in Equation~\eqref{eqn:fully_discrete} is performed in a matrix-free fashion in all numerical examples.

Notably, for the inherited method of Algorithm~\ref{alg:agglomeration_inherited}, the matrix $3\chi_m C_m \frac{M}{2\Delta t} + A$ is assembled only once to generate the coarse-level operators and is then discarded before the time stepping procedure begins. When AMG or Block-Jacobi are employed as preconditioners, the
system matrix must be assembled and handed to the setup phase to construct either the multilevel hierarchy (AMG) or to build a local solver for each block (Block-Jacobi). All the numerical experiments presented in this work directly leverage the matrix-free infrastructure of
the \textsc{deal.II} library for the efficient discontinuous Galerkin operator evaluation~\cite{MatrixFreeDG}.

\section{Numerical results}\label{sec:numerical_experiments_monodomain}

\noindent In this section, we investigate the capabilities of the agglomeration-based preconditioner introduced in Section~\ref{sec:prec_monodomain}
for different scenarios by varying dimensions, geometry, ionic models, and polynomial degrees. As with any multilevel method, our multigrid preconditioner is built upon a few key components:
a hierarchy of meshes, intergrid transfer operators between levels, a smoother, a sequence of operators $\{A_l\}_{l=0}^{L-1}$, and a coarse-grid solver.
The nestedness of the partitions imply that transfer operators are the natural embedding for prolongation and its adjoint for restriction, while the operators on each level are built through the inherited approach described in Section~\ref{sec:prec_monodomain}. The smoother
chosen is a Chebyshev-accelerated Jacobi scheme~\cite{ADAMS2003593} of degree $3$, which uses the inverse of the matrix diagonal on each level (precomputed before solving) and the level operator application for computing residuals in the
Jacobi-type iteration. Optimal parameters are determined by an eigenvalue estimation based on Lanczos iterations. In our simulations, we use $3$ pre- and post-smoothing steps. We compare our multilevel approach against the \textsc{TrilinosML} implementation of AMG~\cite{TrilinosML} and, in order to better position the proposed multigrid method, we also include in
our comparison a Block-Jacobi preconditioner with ILU(0) solver on each block from \textsc{Trilinos}. In this work we use multilevel methods
as preconditioners, since it is well known that multigrid methods are more robust when used as preconditioners rather than solvers~\cite{HariStadlerBirosHighOrderMG}. The preconditioned conjugate-gradient solver applies a single V-cycle (see e.g.~\cite{Trottenberg2001}) of either the agglomerated multigrid or the algebraic multigrid preconditioner. The AMG preconditioner parameters, reported in Table~\ref{tab:amg_params}, are set according
to best practices. All preconditioners are initialized, once and for all, before the time stepping begins.

For all tests, we terminate the PCG iterations when the absolute residual norm falls below $10^{-14}$. At time step $t^{n+1}$, we use the solution from the previous time step $t^n$ as initial guess for the PCG solver. Given the small time step
size $\Delta t$ used in our simulations, this initial guess already closely approximates the solution, requiring only few iterations to meet the convergence criterion. In the left ventricle simulations, we have verified that
the norm of the starting residual of the linear system is approximately $10^{-6}$ throughout the entire time stepping process. Consequently, a stopping tolerance of $10^{-14}$ corresponds to a
eight-order-of-magnitude reduction in the residual norm.

All the numerical experiments are performed using the \textsc{C++} software project \textsc{polyDEAL}~\cite{FEDER2025113773,polydeal}, which is built on top of the finite element library~\textsc{deal.II}~\cite{dealII95,dealIIdesign}. It provides
building blocks for solving PDEs using high-order discontinuous Galerkin methods on polytopic meshes. In particular, it supports the R-tree based agglomeration strategy
and the agglomerated multigrid framework described in the previous section. It is memory distributed using the Message Passing Interface (MPI) standard for parallel computing, and exploits \textsc{p4est}~\cite{p4est} and \textsc{parmetis}~\cite{parmetis} for
distributed mesh refinement and partitioning, while relying on \textsc{Trilinos}~\cite{Trilinos} as main linear algebra backend. All our experiments are publicly available
at the GitHub page of \textsc{polyDEAL}~\cite{polydeal}, and have been run in parallel using different numbers of MPI processes. Detailed instructions to run example programs
are provided in the repository.

\subsection{Two-dimensional test case}
\begin{figure*}[t!]
    \centering
    \begin{minipage}[t]{0.3\linewidth}
        \centering
        \includegraphics[width=\textwidth]{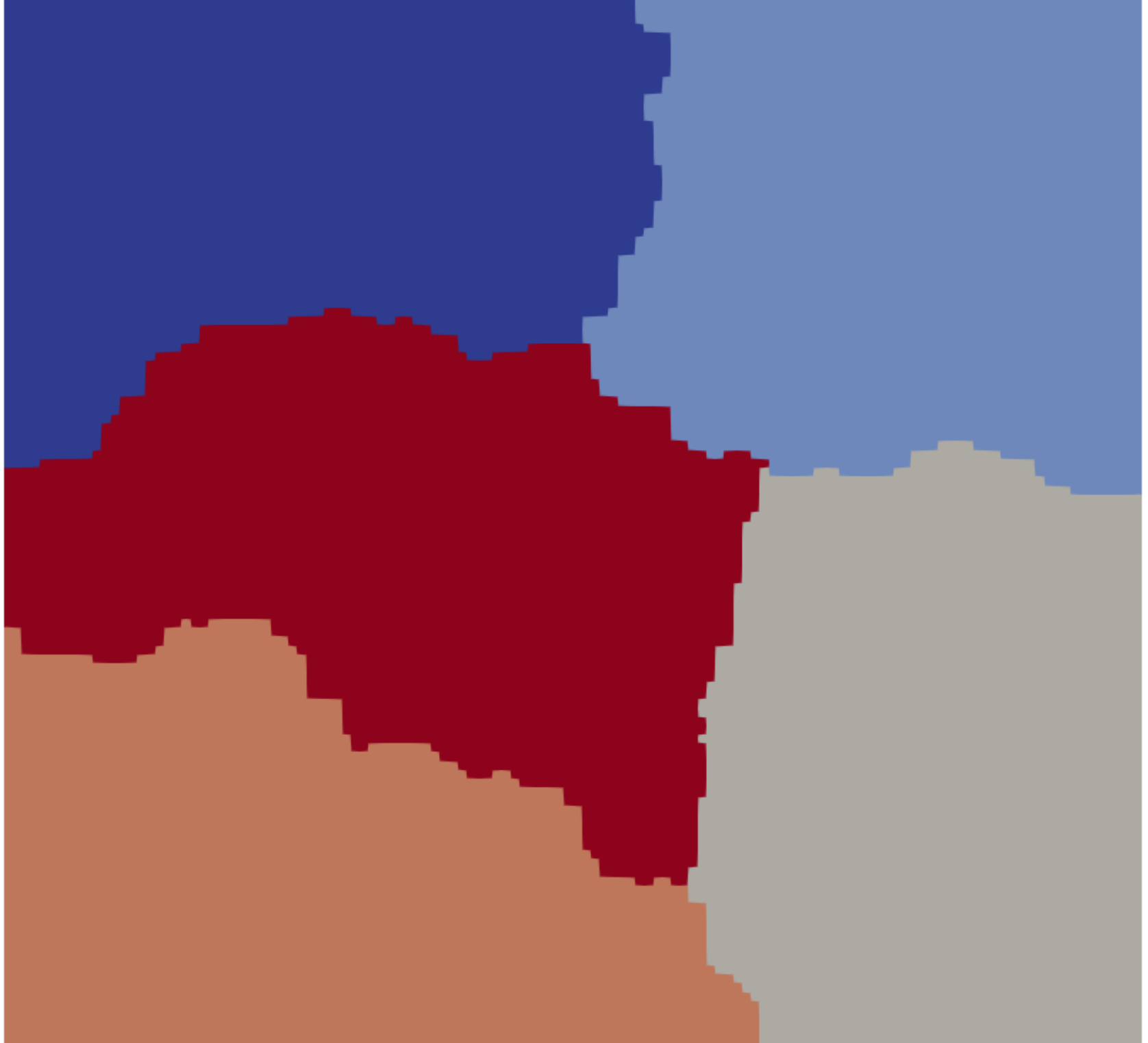}
        \subcaption{MPI partitioning of $\mathcal{T}_h$}
        \label{subfig:2D_partition}

        \centering
        \includegraphics[width=1.\textwidth]{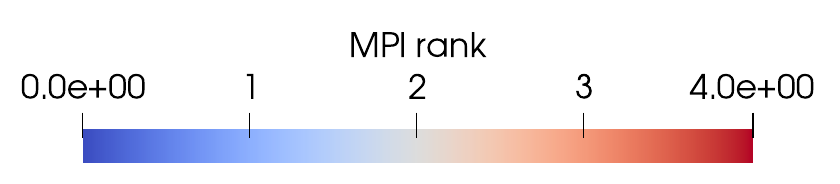}
    \end{minipage}\hfill
    \begin{minipage}[t]{0.63\linewidth}
        \centering
        \begin{minipage}[t]{0.44\linewidth}
            \centering
            \includegraphics[width=\textwidth]{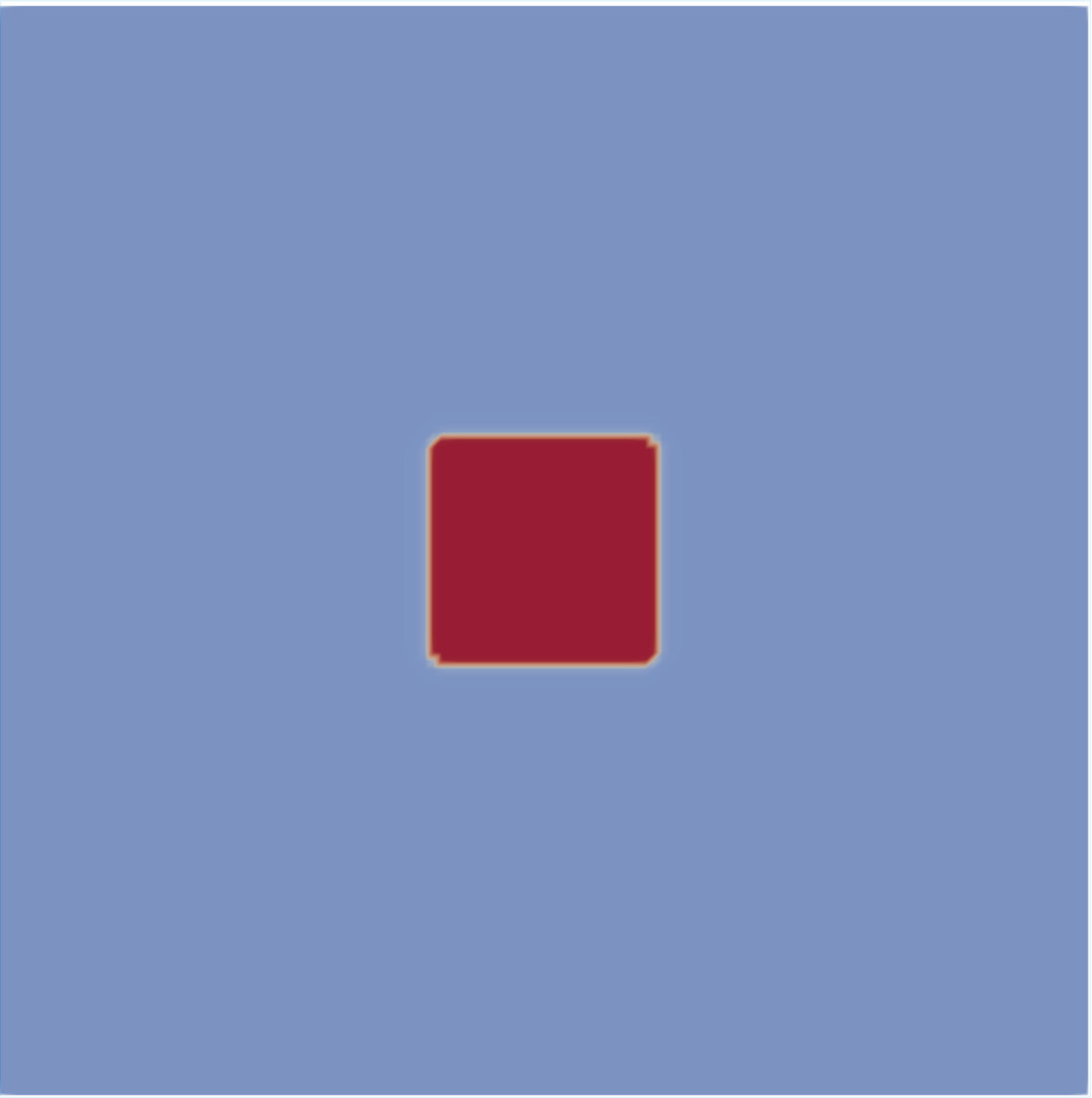}
            \subcaption{$t=0.04$s}
            \label{subfig:2D_t004}
        \end{minipage}\hfill
        \begin{minipage}[t]{0.44\linewidth}
            \centering
            \includegraphics[width=\textwidth]{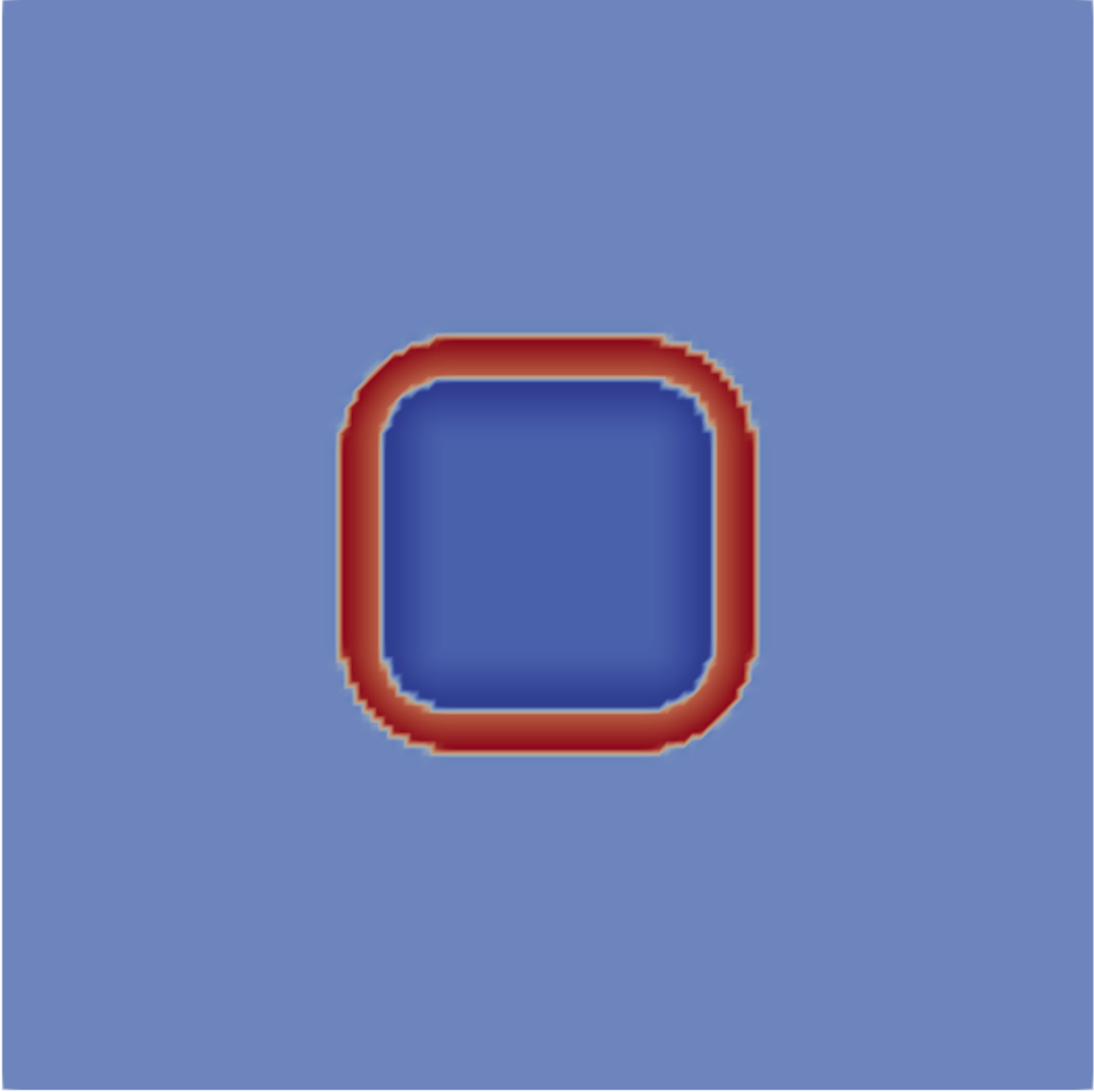}
            \subcaption{$t=0.16$s}
            \label{subfig:2D_t016}
        \end{minipage}

        \centering
        \includegraphics[width=0.5\textwidth]{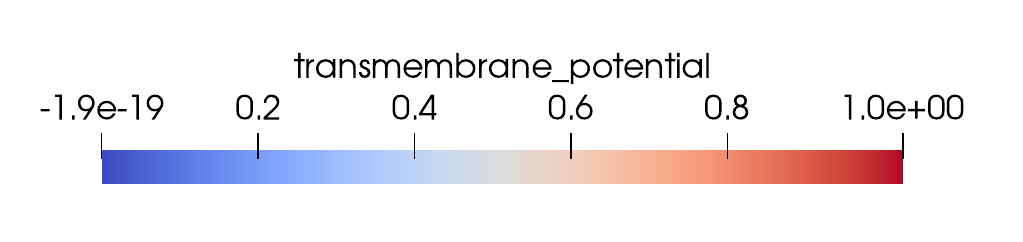}
    \end{minipage}

    \caption{(\subref{subfig:2D_partition}) The computational mesh partitioned using \textsc{parmetis} among five MPI processes. To build the multigrid hierarchy, aggregates are generated using the
    locally owned cells within each MPI rank. Snapshots of the transmembrane potential $\boldsymbol{U}$ at (\subref{subfig:2D_t004}) $t=0.04$s and (\subref{subfig:2D_t016}) $t=0.16$s after the external application of the current $\mathcal{I}_{\text{app}}(\boldsymbol{x},t)$ in Equation~\eqref{eqn:I_app_2D}.}
    \label{fig:2D_results}
\end{figure*}

We start with the following two-dimensional test, taken from~\cite{botta2024highorderdiscontinuousgalerkinmethods}, aiming at assessing both the validity of our
implementation with a test case which exhibits both the depolarization and repolarization phases, as well as to test the proposed preconditioning strategy. We consider
Problem~\eqref{eqn:monodomain}, over $\Omega \!=\! (0,1)^2$, coupled with the FitzHugh-Nagumo ionic model described by the ODE~\eqref{eqn:FitzHugNagumo}. The parameters of the
ionic model are identical to the ones used in~\cite{botta2024highorderdiscontinuousgalerkinmethods}, which we report in Table~\ref{tab:FHZ_params}
for the sake of completeness. The unit square is refined $7$ times, resulting in a fine mesh $\mathcal{T}_h$ consisting of $16,384$ quadrilateral elements. The number of degrees of freedom for different polynomial degrees is reported in Table~\ref{tab:dofs_2D}. The mesh $\mathcal{T}_h$ is then partitioned through the mesh partitioner \textsc{parmetis}~\cite{parmetis}, where each MPI process
owns a subdomain made of contiguous elements. This implies that, in general, the resulting MPI partitions are not describing structured grids anymore, as shown in Fig.~\ref{subfig:2D_partition}. As previously mentioned in Section~\ref{subsec:rtree-agglomeration}, we generate
agglomerates \emph{within} each locally owned partition. As such, this test is representative of more general configurations in which unstructured geometries are considered. Indeed, since
the original mesh is structured, the grids generated by the R-tree algorithm without partitioning
the mesh among processes\footnote{Which means a serial run with $1$ MPI process.} would be a structured hierarchy of coarser meshes.

The applied current is set to be: \begin{equation}\label{eqn:I_app_2D}\mathcal{I}_{\text{app}}(\boldsymbol{x},t) = 2 \cdot 10^6\mathbbm{1}_{[0.4,0.6]}(x)\mathbbm{1}_{[0.4,0.6]}(y)\mathbbm{1}_{[0,10^{-3}]}(t),\end{equation}
where $\mathbbm{1}_{[a,b]}(\cdot)$ is the indicator function over the interval $[a,b]$. Such definition corresponds to a temporary
electric shock localized in the square $[0.4,0.6]^2$. We solve the problem until time $T=0.4\si{\second}$, with step size $\Delta t=10^{-4}\si{\second}$. The potential $u$ and the gating variable $w$ are set at rest.\footnote{Using the FitzHug-Nagumo model, $\boldsymbol{w}$ is a scalar unknown.} We show the activation and the subsequent
propagation of the transmembrane potential $u$ in Figs.~\ref{subfig:2D_t004} and~\ref{subfig:2D_t016}, obtained with polynomial degree $p\!=\!2$. As shown, the wavefront is accurately captured by the DG discretization.

\begin{figure}[!t]
    \centering

    \includegraphics[width=\textwidth]{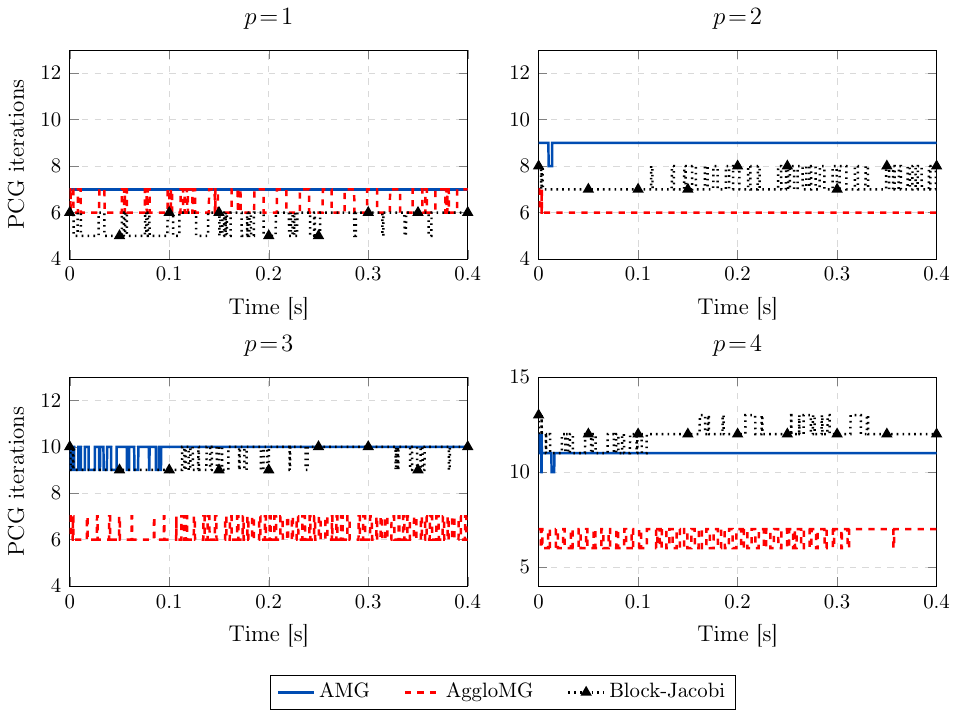}

    \caption{Number of PCG iterations per time step for Problem~\eqref{eqn:monodomain}, comparing AMG, agglomerated multigrid (AggloMG), and Block-Jacobi for polynomial degrees $p=1,2,3,4$.}
    \label{fig:iterates_2D_monodomain}
\end{figure}

As coarse grid solver, we use the parallel sparse direct solver \textsc{mumps}~\cite{mumps}. The settings for the agglomerated multigrid preconditioner are shown in Table~\ref{tab:agglo_mg_config_2d}. We report in
Fig.~\ref{fig:iterates_2D_monodomain} the number of outer iterations required by conjugate gradient for
each time step, varying the polynomial degree $p$ from $1$ to $4$. For all polynomial degrees, the agglomeration-based multilevel strategy achieves lower
iterations counts compared to the AMG preconditioner, even for the lowest-order case $p\!=\!1$, which is the
most favorable for AMG solvers. For $p\!=\!4$, we set the AMG aggregation threshold parameter to $0.1$, as the default value of $0.2$ failed to build a multigrid hierarchy. The Block-Jacobi preconditioner shows iteration counts similar to AMG, with an increase with respect to the polynomial degree $p$.

Overall, the agglomeration-based preconditioner displays a clear robustness with respect to the polynomial degree $p$, while the AMG and Block-Jacobi iterations are increasingly higher as $p$ grows, doubling the
number of iterations of the agglomeration-based preconditioner when $p\!=\!4$. Conversely, the agglomeration-based
preconditioner yields for every $p$ a number of iterations that is essentially constant and ranging from $6$ to $7$, confirming the effectiveness of the
proposed strategy.

\begin{table}[!t]
    \centering
    \begin{tabular}{@{}cc@{}}
        \toprule
        Polynomial degree $p$ & DoFs    \\
        \midrule
        $1$                   & 65,536  \\
        $2$                   & 147,456 \\
        $3$                   & 262,144 \\
        $4$                   & 409,600 \\
        \bottomrule
    \end{tabular}
    \caption{Number of degrees of freedom for the two-dimensional test with polynomial degrees $p\!=\!1,\ldots,4$.}
    \label{tab:dofs_2D}
\end{table}

\subsection{Idealized three-dimensional test case}
\begin{figure}[!t]
    \centering
    \begin{minipage}[c]{.42\linewidth}
        \centering
        \includegraphics[width=0.7\textwidth]{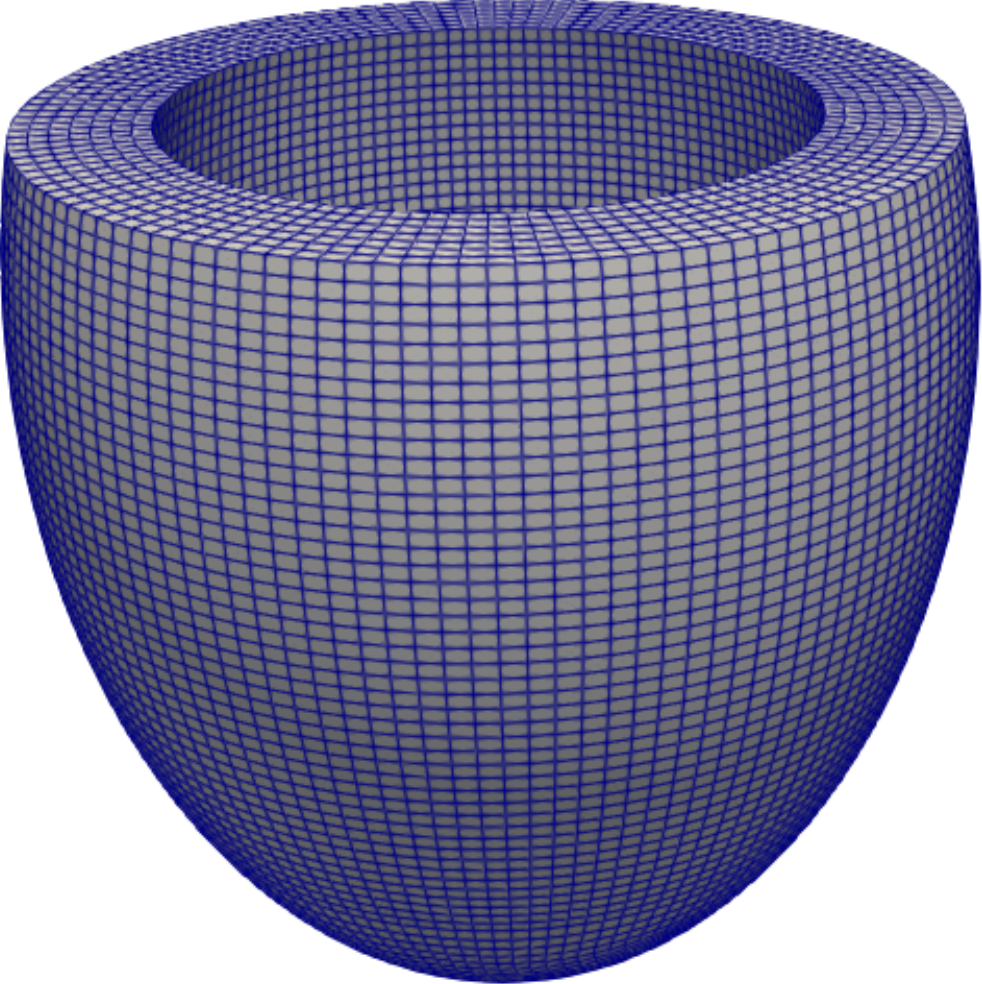}
        \subcaption{Idealized mesh of the left ventricle.}
        \label{fig:left_ventricle_ideal}
    \end{minipage}%
    \hfill
    \begin{minipage}[c]{.5\linewidth}
        \centering

        \includegraphics[width=\textwidth]{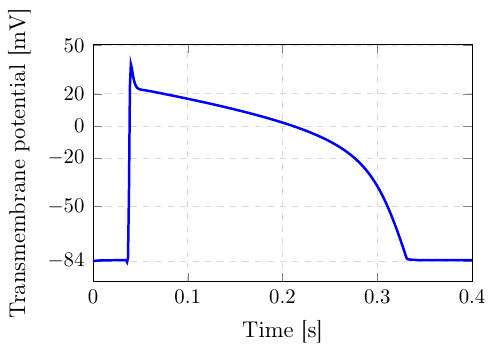}

        \subcaption{Evolution of the transmembrane potential at a point.}
        \label{fig:potential_evolution}
    \end{minipage}
    \caption{Idealized left ventricle test case: (\subref{fig:left_ventricle_ideal}) hexahedral mesh of the ellipsoid and (\subref{fig:potential_evolution}) time evolution of the transmembrane potential at a selected point.}
    \label{fig:idealized_ventricle_combined}
\end{figure}

\begin{figure}[h!t]
    \centering
    \begin{subfigure}{0.3\linewidth}
        \centering
        \includegraphics[width=\linewidth]{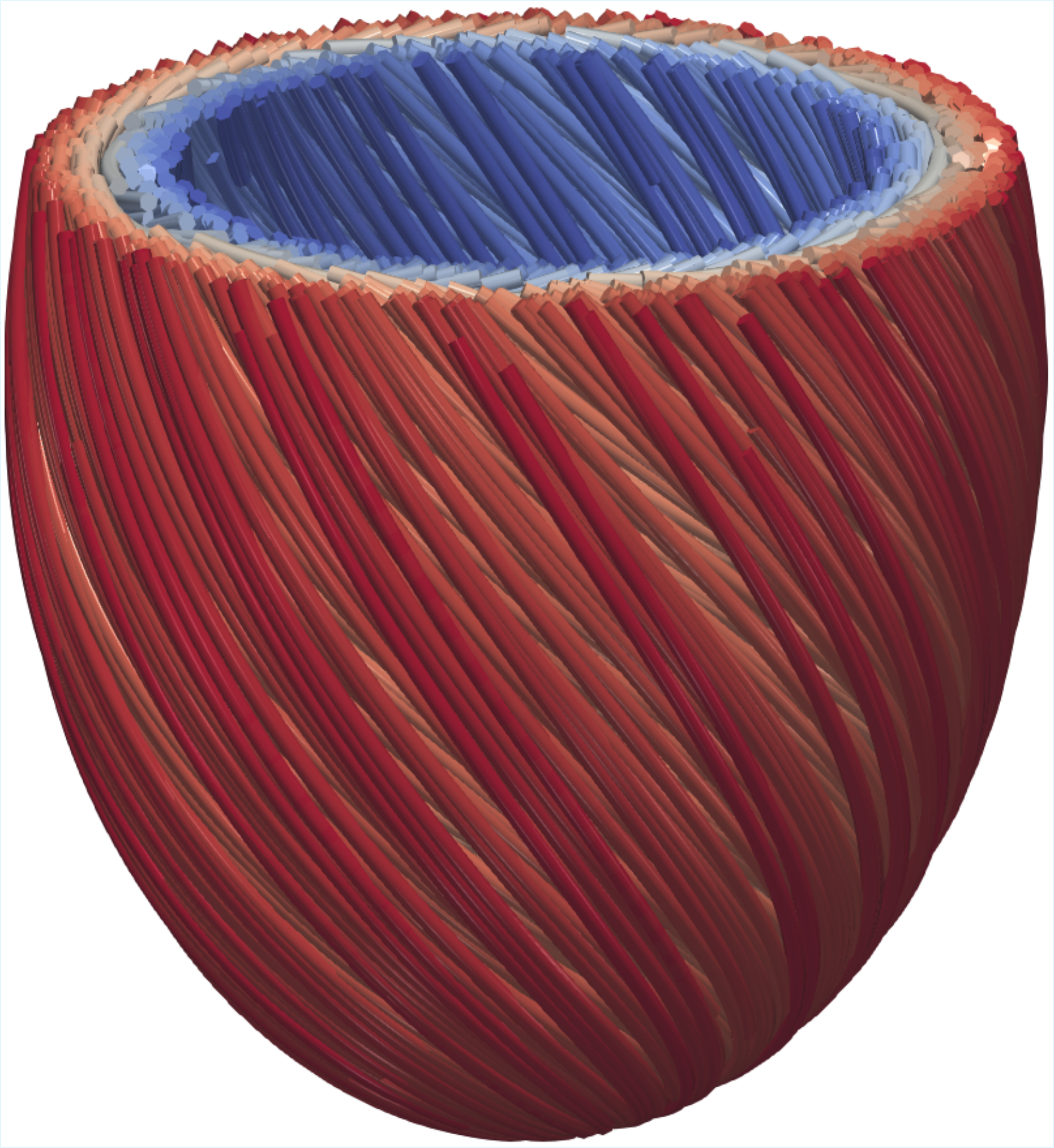}
        \caption{}
        \label{subfig:fibers_idealized}
    \end{subfigure}
    \hspace{3cm}
    \begin{subfigure}{0.3\linewidth}
        \centering
        \includegraphics[width=\linewidth]{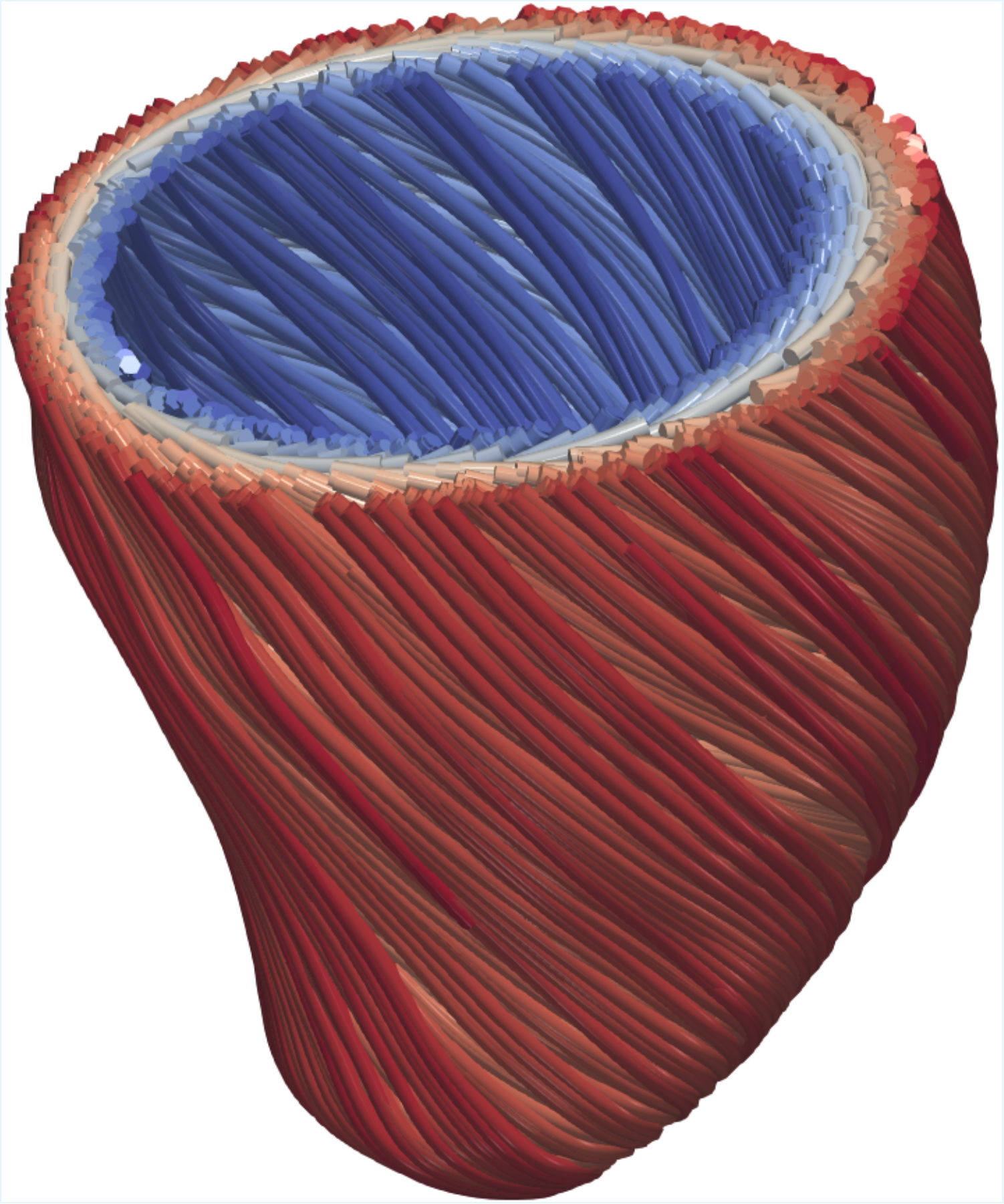}
        \caption{}
        \label{subfig:fibers_realistic}
    \end{subfigure}

    \centering
    \includegraphics[width=0.4\textwidth]{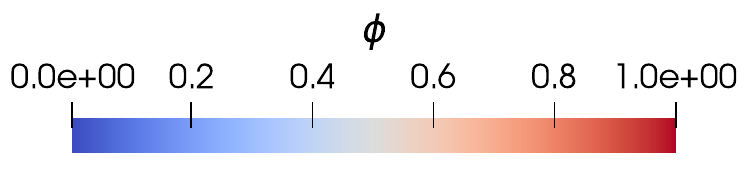}

    \caption{Fiber field computed using the Laplace-Dirichlet Rule-Based Methods~\cite{Africa2023lifexfiber,PIERSANTI2021113468} visualized as streamlines: (\subref{subfig:fibers_idealized}) Idealized left ventricle; (\subref{subfig:fibers_realistic}) Realistic left ventricle. The solution
        $\phi$ is the transmural function, where $\phi\!=\!0$ at the endocardium and $\phi\!=\!1$ at the epicardium.}
    \label{fig:fibers}
\end{figure}

As outlined in Section~\ref{sec:prec_monodomain}, the main drawback of building the coarse operators through an inherited approach is that the
number of iterations depends on the number of levels in the multigrid hierarchy~\cite{AntoniettiSartiVeraniDGhpMG}. For this reason, we investigate this
dependence for the linear system of equations arising from Equation~\eqref{eqn:fully_discrete}.
To provide a meaningful test, we consider the idealized left ventricle shown in Fig.~\ref{fig:left_ventricle_ideal}, consisting of 407,904 hexahedral elements. The fiber field computed on this geometry is shown in Fig.~\ref{subfig:fibers_idealized}. Partitioning the mesh $\mathcal{T}_h$ among
$5$ processes, we obtain the levels with the cardinalities reported in Table~\ref{tab:left_ventricle_ideal}. We solve the monodomain problem
with different number of levels and polynomial degrees $p$ from 1 to 3 in order to
verify whether the number of PCG iterations is affected by the number of levels in the hierarchy. We use a time step $\Delta t \!=\!10^{-4}\si{\second}$, solving Problem~\eqref{eqn:monodomain} until the final time
$T\!=\!0.4\si{\second}$. The settings for the agglomerated multigrid preconditioner are shown in Table~\ref{tab:agglo_mg_config_3d}, while we use the
Bueno-Orovio ionic model defined by the system of ODEs in~\eqref{eqn:BuenoOrovio}, with parameters reported in Table~\ref{tab:BO_params}.
The external stimulus $\mathcal{I}_{\text{app}}(\boldsymbol{x},t)$ is localized at three points of the ventricle for $3\si{\milli\second}$ with a
magnitude of $300\si{\per\second}$. As initial data we take $u\!=\!0$ for the transmembrane
potential, while the initial value for the gating variables are $w_0\!=\!w_1\!=\!1$, and $w_2\!=\!0$. The transmembrane potential in millivolts is obtained by post-processing
the numerical solution using the formula $u_{\si{\milli\volt}}\!=\!(85.7u - 84) [\si{\milli\volt}]$~\cite{BUENOOROVIO2008544}. We show in Fig.~\ref{fig:potential_evolution} the time
evolution of the transmembrane potential at a selected point of the ventricle, where the typical depolarization and repolarization phases are clearly visible.

We report in Table~\ref{tab:iterations_level_statistics} the statistics of the number of PCG iterations per time step for all the considered configurations in terms of polynomial degrees and number of levels. The mean
values and the standard deviations confirm that for this problem the number of iterations is essentially independent of the number of levels employed. Notably, iteration counts are
low and robust for all polynomial degrees.

This experimental observations are not in contradiction with the theoretical results in~\cite{AntoniettiSartiVeraniDGhpMG,BOTTI2017382}, as here
we are not solving a classical Poisson system, but rather the monodomain problem, where the presence
of the dominant term $\frac{M}{\Delta t}$ shifts the spectrum of the system away from zero. We remark that the choice of $\Delta t$ is dictated by the need to accurately capture the
upstroke across the depolarization phase of the action potential, which typically requires small time steps~\cite{AFRICA2023111984,FranzonePavarinoParallelMonodomain}. In particular, $\Delta t$ is bounded from above by the stimulus duration ($3\si{\milli\second}$ in our case). We have further
verified (Table~\ref{tab:multigrid_iterations_delta_t}) that the number of PCG iterations increases with $\Delta t$ for all the solvers considered (consistent with the improvement in conditioning as $\Delta t$ decreases), the increase being most pronounced for Block-Jacobi. Across the admissible range
of time steps and for all polynomial degrees, the agglomeration-based multigrid preconditioner remains robust
and consistently yields the lowest iteration counts.

In the forthcoming section, we provide a comparison with AMG and Block-Jacobi on a three-dimensional test case with a realistic geometry.

\begin{table}[!t]
    \centering
    \begin{minipage}[t]{0.35\textwidth}
        \centering
        \begin{tabular}{@{}lr@{}}
            \toprule
            \multicolumn{2}{c}{Idealized ventricle}        \\
            \multicolumn{2}{c}{($5$ processes)}            \\
            \midrule
            Level index $l$ & $\text{Card}(\mathcal{T}_l)$ \\
            \midrule
            $l=0$           & 407,904                      \\
            $l=1$           & 50,990                       \\
            $l=2$           & 6,375                        \\
            $l=3$           & 800                          \\
            $l=4$           & 100                          \\
            \bottomrule
        \end{tabular}
        \caption{Coarsened hierarchy $\{\mathcal{T}_l\}_{l=0}^4$ for the mesh representing the left ventricle, partitioning the original mesh across $5$ processes.}
        \label{tab:left_ventricle_ideal}
    \end{minipage}%
    \hfill
    \begin{minipage}[t]{0.58\textwidth}
        \centering
        \begin{tabular}{@{}ccccccc@{}}
            \toprule
            \multirow{2}{*}{$p$} & \multirow{2}{*}{DoFs}       & \multirow{2}{*}{\shortstack{No. of                              \\Levels ($L$)}} & \multicolumn{4}{c}{PCG Iterations} \\
            \cmidrule(lr){4-7}
                                 &                             &                                    & Mean & Std Dev & Min & Max \\
            \midrule
            \multirow{3}{*}{$1$}
                                 & \multirow{3}{*}{3,263,232}  & 2                                  & 7.6  & 0.6     & 6   & 9   \\
                                 &                             & 3                                  & 7.5  & 0.6     & 6   & 9   \\
                                 &                             & 4                                  & 7.3  & 0.9     & 5   & 9   \\
            \midrule
            \multirow{3}{*}{$2$}
                                 & \multirow{3}{*}{11,013,408} & 2                                  & 8.1  & 0.7     & 7   & 10  \\
                                 &                             & 3                                  & 8.1  & 0.7     & 7   & 10  \\
                                 &                             & 4                                  & 8.1  & 0.8     & 7   & 10  \\
            \midrule
            \multirow{3}{*}{$3$}
                                 & \multirow{3}{*}{26,105,856} & 2                                  & 8.7  & 0.9     & 7   & 11  \\
                                 &                             & 3                                  & 8.5  & 1.0     & 7   & 11  \\
                                 &                             & 4                                  & 8.5  & 1.0     & 7   & 11  \\
            \bottomrule
        \end{tabular}
        \caption{Statistics (mean value, standard deviation, and min-max) of PCG iterations per time step for the idealized ventricle test case across different numbers of multigrid levels $L$ and polynomial degrees $p$.}
        \label{tab:iterations_level_statistics}
    \end{minipage}
\end{table}

\subsection{Realistic three-dimensional test case}
In the following experiment, we use the realistic CAD-modeled left-ventricle mesh $\mathcal{T}_h$ introduced in Section~\ref{subsec:rtree-agglomeration}, together with the same model and solver parameters as in the previous idealized
test case, except for the conductivities which are set to $\sigma_l = 10\times10^{-5}$ [\si{\meter\squared\per\second\squared}], $\sigma_t = 2\times10^{-5}$ [\si{\meter\squared\per\second\squared}], and $\sigma_n = 2\times10^{-5}$ [\si{\meter\squared\per\second\squared}]. The fiber field computed on this mesh is shown in Fig.~\ref{subfig:fibers_realistic}.

\begin{figure}
    \centering
    \begin{subfigure}{0.30\linewidth}
        \centering
        \includegraphics[width=\linewidth]{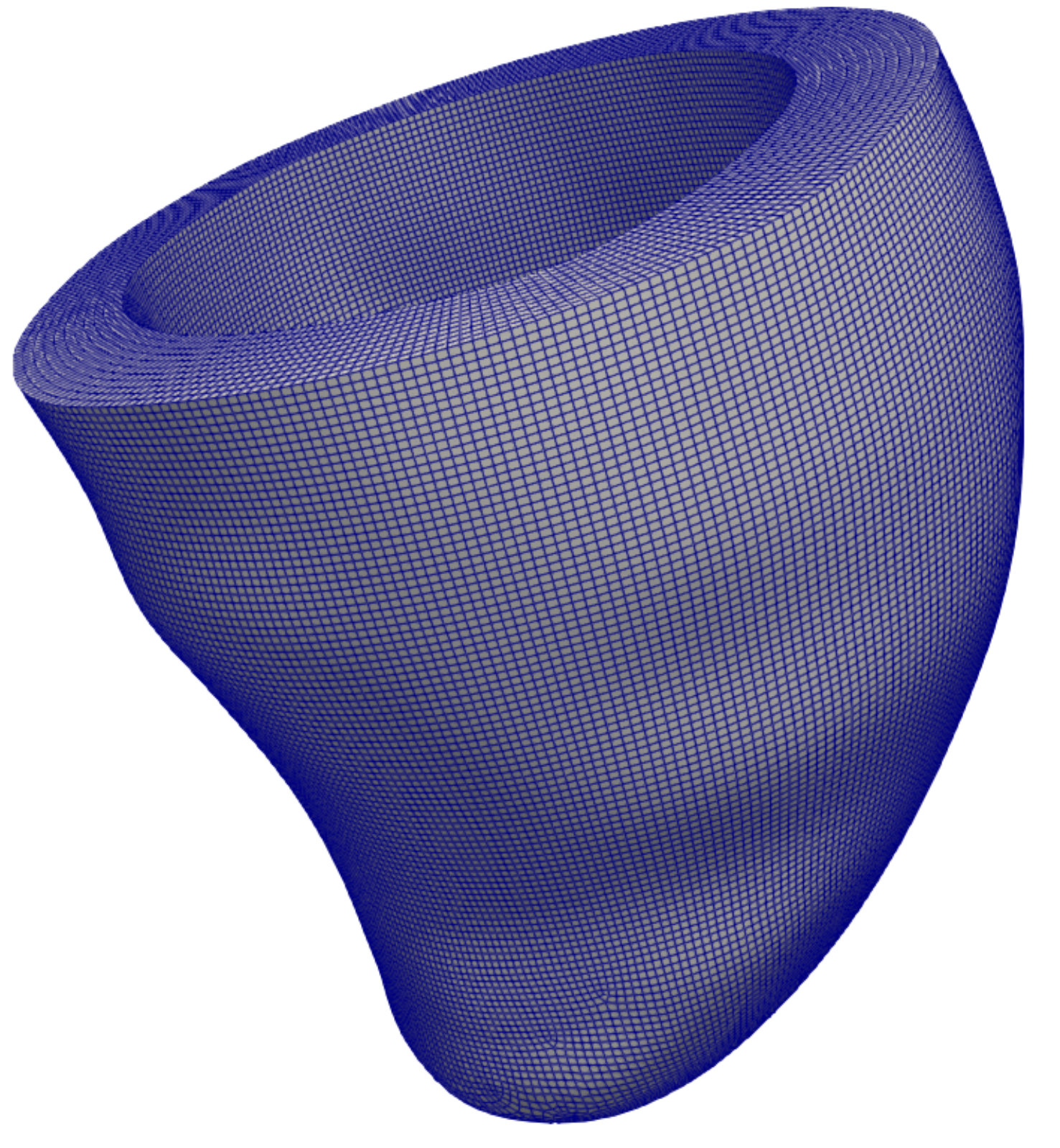}
        \caption{}
        \label{subfig:left_ventricle_mesh}
    \end{subfigure}
    \hspace{4cm}
    \begin{subfigure}{0.30\linewidth}
        \centering
        \includegraphics[width=\linewidth]{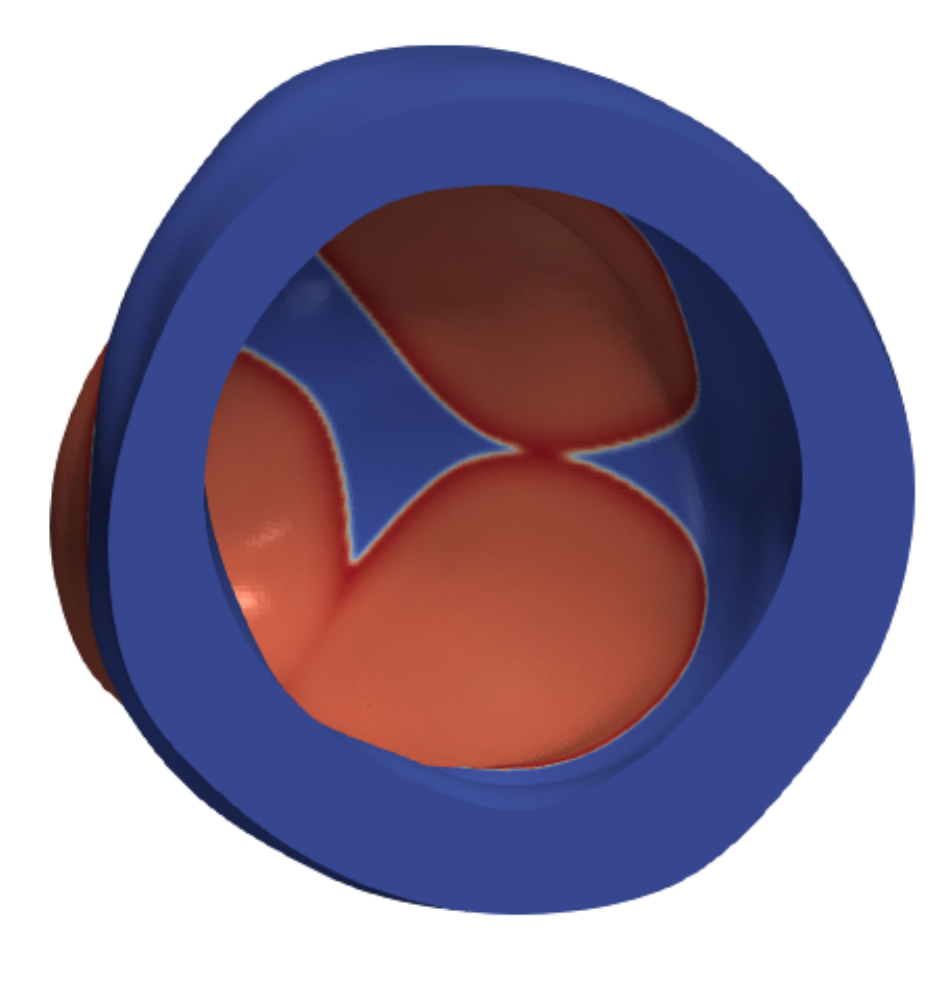}
        \caption{}
        \label{subfig:left_ventricle_impulses}
    \end{subfigure}

    \caption{
        (\subref{subfig:left_ventricle_mesh}) Hexahedral mesh $\mathcal{T}_h$ representing a realistic left ventricle.
        (\subref{subfig:left_ventricle_impulses}) Propagation of the transmembrane potential $\boldsymbol{U}$.
    }
    \label{fig:left_ventricle}
\end{figure}

A visualization of the mesh and of the propagation of the transmembrane potential $\boldsymbol{U}$ is shown in Figure~\ref{fig:left_ventricle}. In our investigation, we consider
polynomial degrees from $p\!=\!1$ to $p\!=\!3$.

The multigrid hierarchy has been generated in parallel using the agglomeration procedure outlined in Algorithm~\ref{alg:r_tree_algo}. The number of mesh elements and degrees of freedom (with $\mathcal{Q}^1$ elements) at
each level $l\in \{0,\ldots,L-1\}$ is shown in Table~\ref{tab:Ratio_meshes_ventricle} when the
full ventricle geometry is distributed among an increasing number of MPI processes.
Notably, the decrease of the global number of elements between consecutive levels has a ratio very close to $8$, which is the value one would obtain by halving the mesh
step size of a uniform hexahedral grid. We observe again good balancing in terms of workload per process and ratios between the cardinality of consecutive grids. This is possible thanks to the combination of the graph partitioner \textsc{parmetis}, applied to
the initial grid $\mathcal{T}_h$, and the balanced number of elements achieved by the agglomeration routine.

To ensure efficient multigrid convergence, the hierarchy must contain enough levels $L$ such that
the coarsest grid has an adequate size, preventing the coarse-level solve from becoming a bottleneck in the V-cycle.
Direct solvers on the coarsest level pose scalability challenges in distributed-memory settings, as the resulting small system size leads to communication overhead that far exceeds computational work.
One approach, implemented in \textsc{Hypre}~\cite{hypre}, restricts the coarse solve to a subset of MPI ranks
while idling the remaining processes, followed by redistribution of the solution.
For our three-dimensional experiments, we instead apply as coarse solver the PCG method with an algebraic multigrid
V-cycle preconditioner, limiting the maximum size of the coarsest-level matrix to a sufficiently large size so that all processes are
involved in the computation. Since PCG is employed as the coarse solver, our preconditioner no longer represents a fixed linear operator, which in principle calls for
the \emph{flexible} variant of the conjugate gradient method\footnote{In practice, the classical conjugate gradient method required the same number of iterations, likely due to the rapid and accurate convergence of the coarse-level PCG iterations.}.

\begin{figure*}[htbp]
    \subfloat[$t=0.04\si{\second}$]{%
        \includegraphics[width=.23\linewidth]{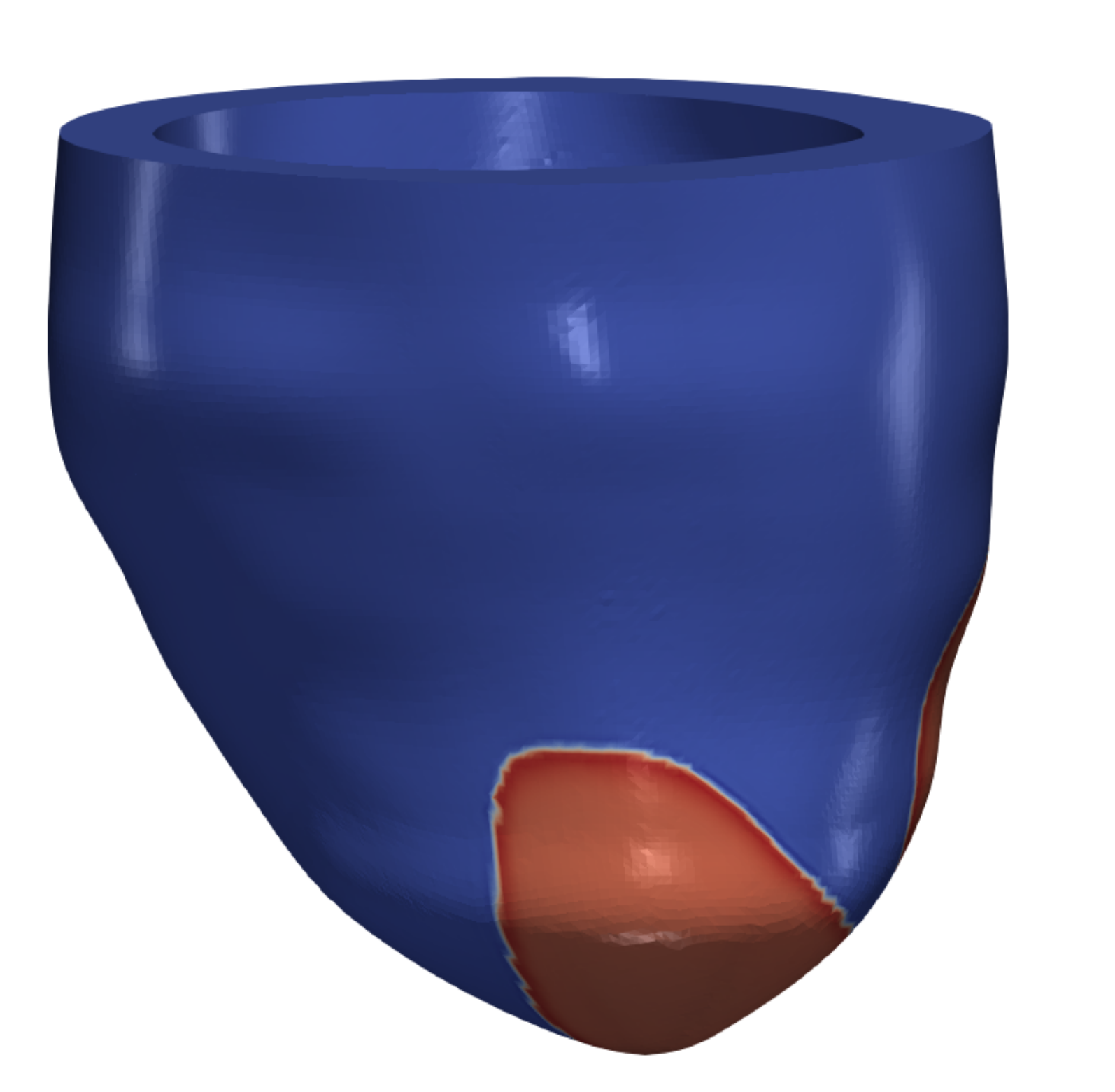}%
        \label{subfig:a}%
    }\>\>\>\>\>\>\>\>\>\>\>\>\>\>\>\>\>\>\>\>\>\>\>\>
    \subfloat[$t=0.08\si{\second}$]{%
        \includegraphics[width=.23\linewidth]{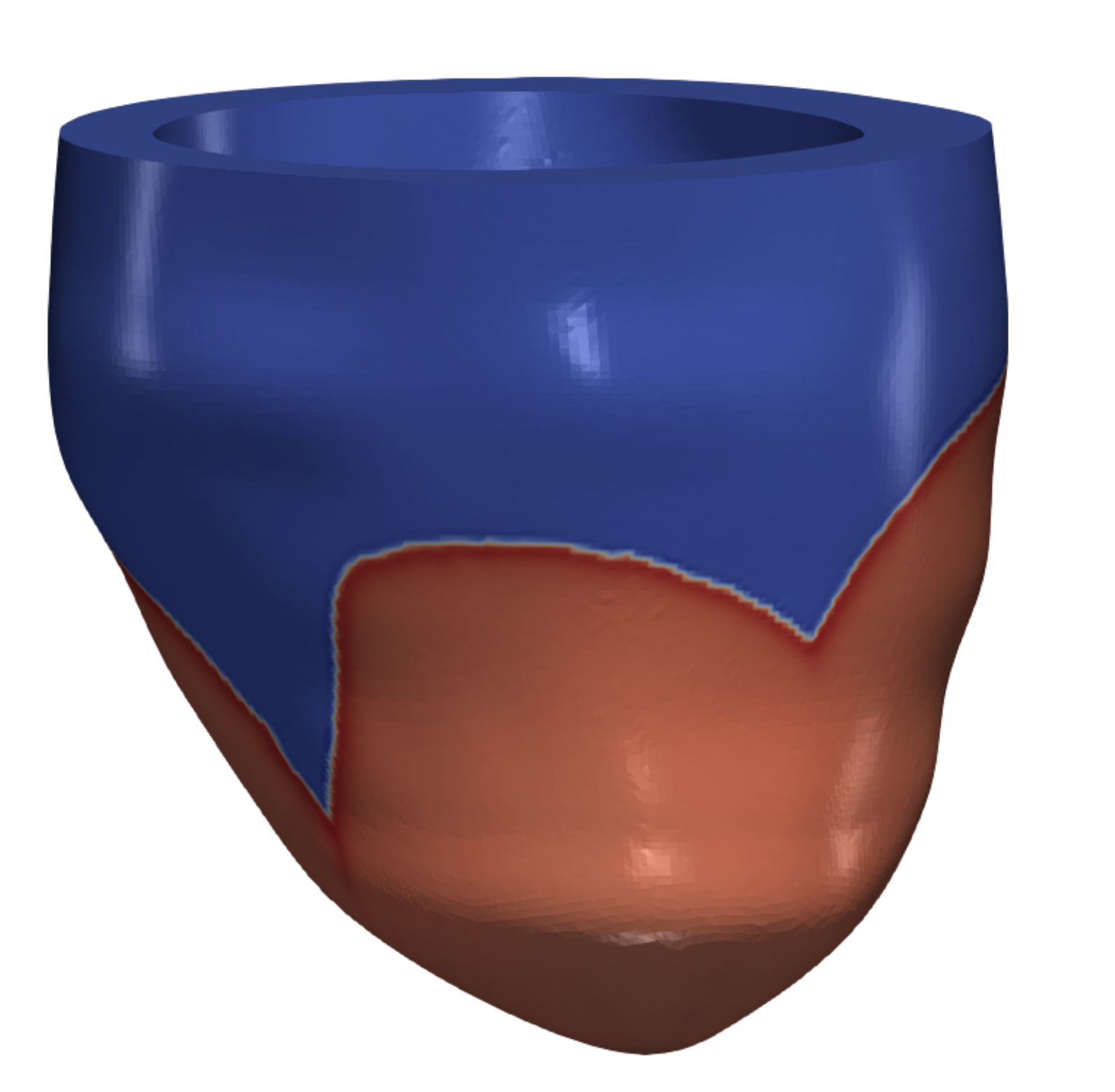}%
        \label{subfig:b}%
    }\>\>\>\>\>\>\>\>\>\>\>\>\>\>\>\>\>\>\>\>\>\>\>\>
    \subfloat[$t=0.2\si{\second}$]{%
        \includegraphics[width=.23\linewidth]{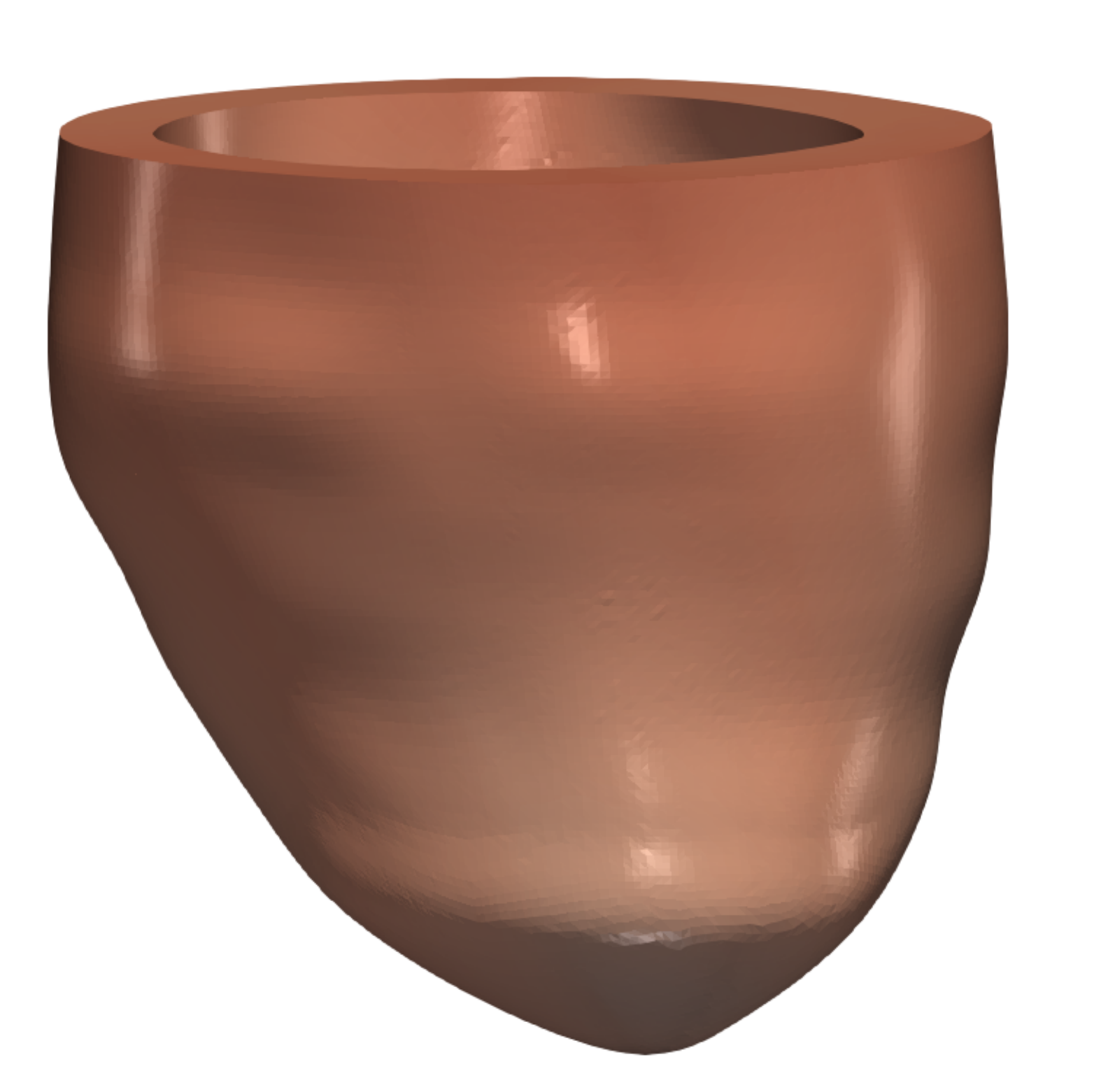}%
        \label{subfig:c}%
    }\\
    \subfloat[$t=0.28\si{\second}$]{%
        \includegraphics[width=.23\linewidth]{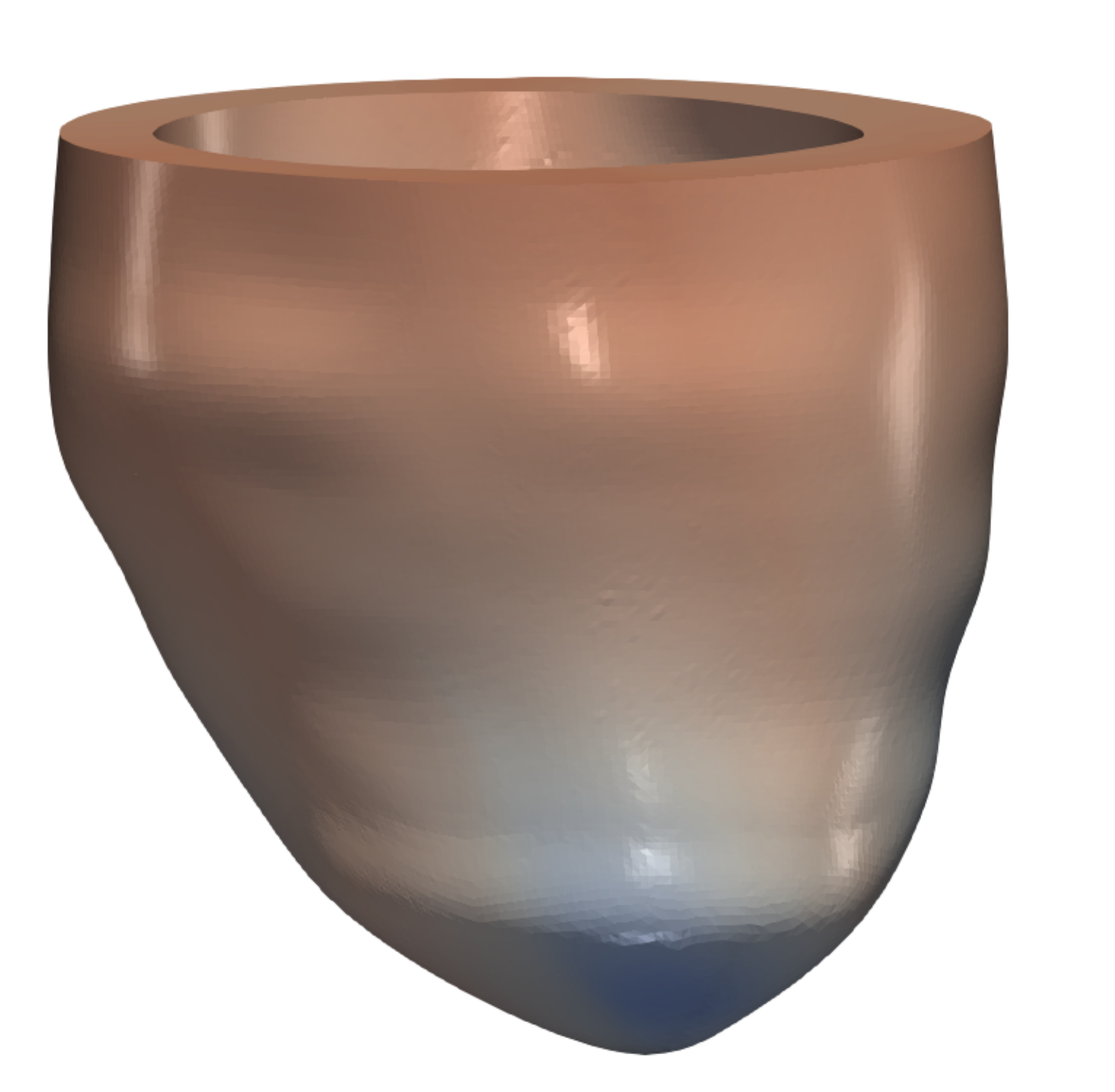}%
        \label{subfig:d}%
    }\>\>\>\>\>\>\>\>\>\>\>\>\>\>\>\>\>\>\>\>\>\>\>\>
    \subfloat[$t=0.32\si{\second}$]{%
        \includegraphics[width=.23\linewidth]{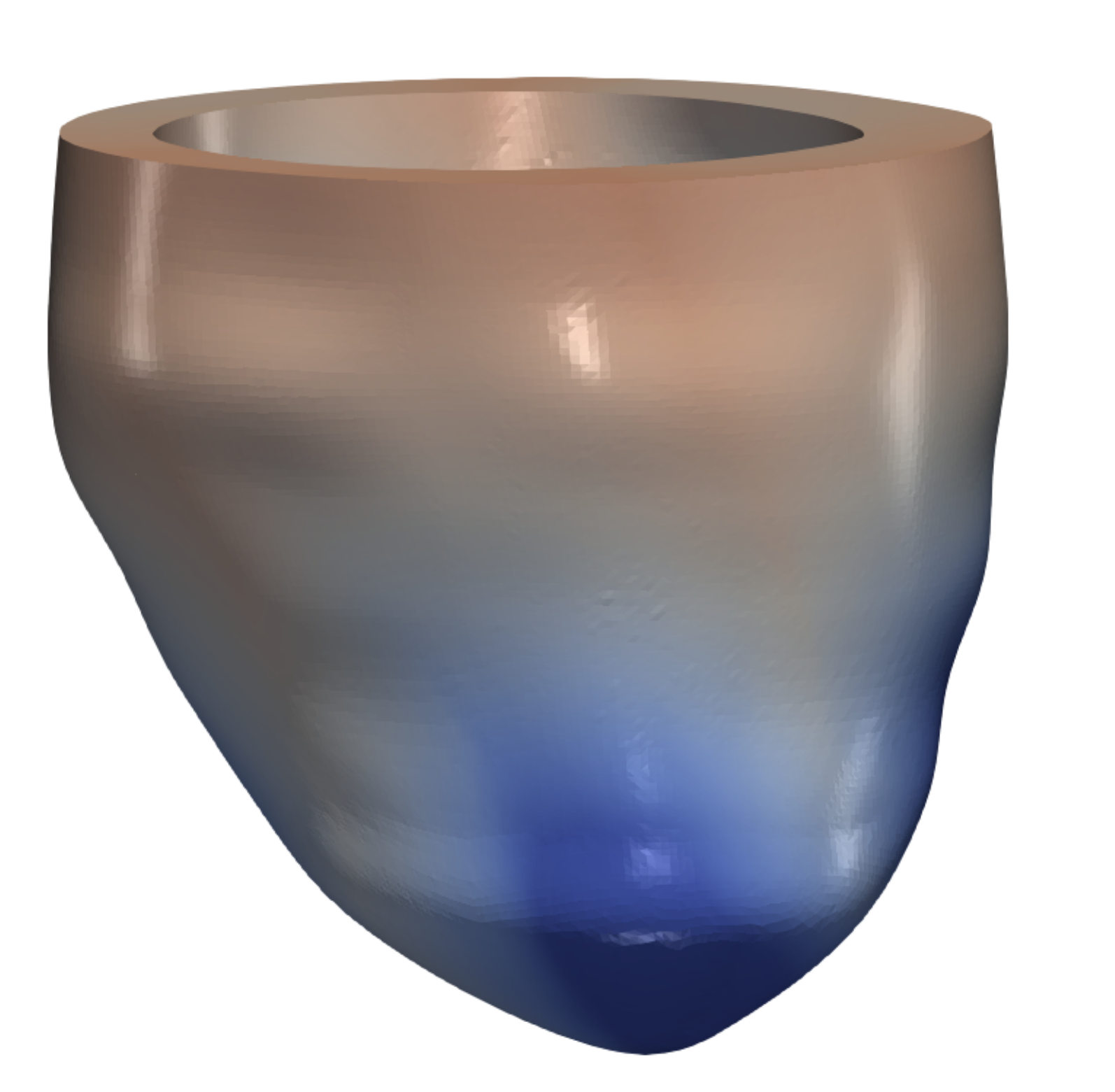}%
        \label{subfig:e}%
    }\>\>\>\>\>\>\>\>\>\>\>\>\>\>\>\>\>\>\>\>\>\>\>\>
    \subfloat[$t=0.36\si{\second}$]{%
        \includegraphics[width=.23\linewidth]{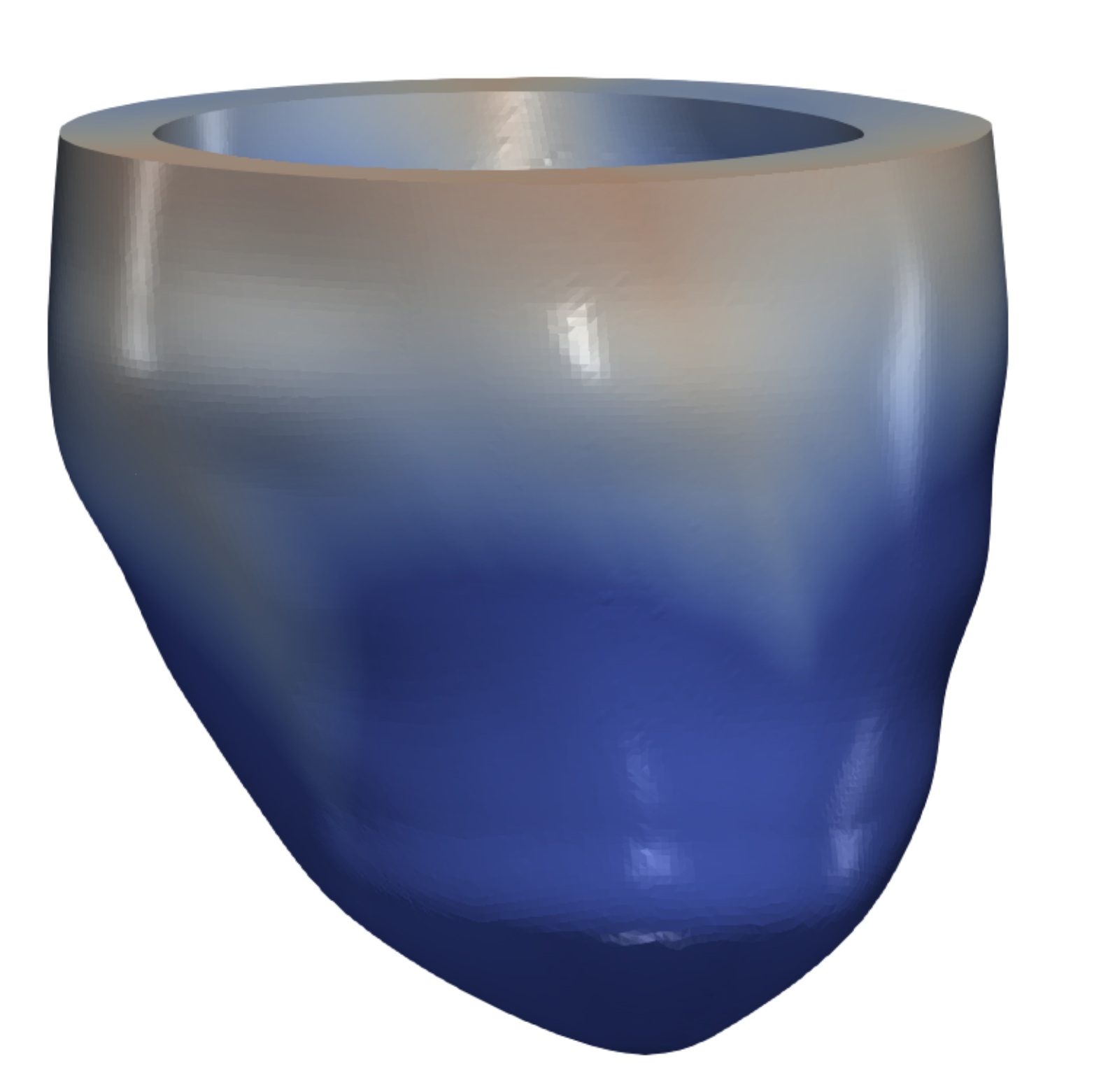}%
        \label{subfig:f}%
    }\\
    \centering
    \includegraphics[width=0.45\textwidth]{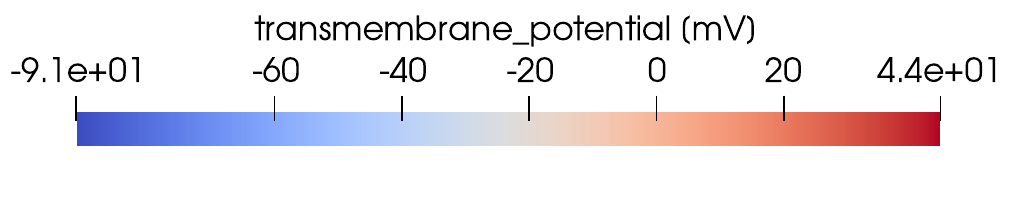}%
    \caption{Snapshots of the transmembrane potential $\boldsymbol{U}$ at selected time steps after the external application of the current $\mathcal{I}_{\text{app}}(\boldsymbol{x},t)$ at three locations inside the ventricle.}
    \label{fig:3D_transmembrane}
\end{figure*}

\begin{figure}[!t]
    \centering

    \includegraphics[width=\textwidth]{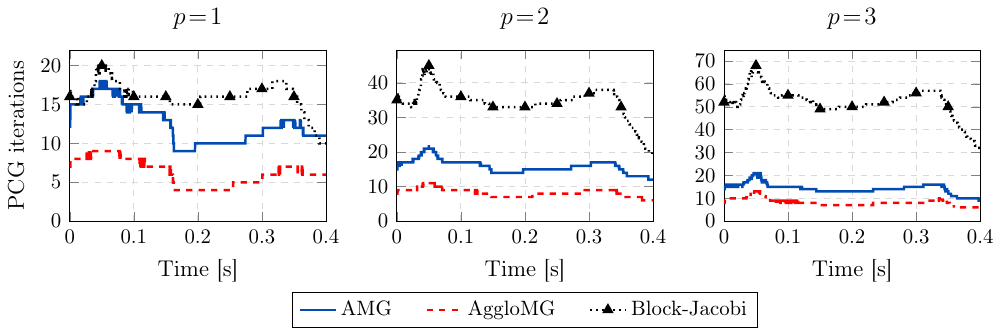}

    \caption{Number of PCG iterations per time step for Problem~\eqref{eqn:monodomain}, comparing AMG, agglomerated multigrid (AggloMG), and Block-Jacobi for polynomial degrees $p\!=\!1,2,3$ for the three-dimensional
        ventricle test.}
    \label{fig:iterates_3D_monodomain}
\end{figure}

We show some snapshots of the potential at selected time steps in Fig.~\ref{fig:3D_transmembrane}, where the propagation of the steep wavefront is visible. We report the iterations' path per time step with linear, quadratic, and cubic polynomial degrees in Fig.~\ref{fig:iterates_3D_monodomain}. For all polynomial degrees and
preconditioners, the number of PCG iterations have a maximum during the initial propagation of the potential, and a successive
drop around $t\!=\!0.12\si{\second}$, when the solution does not exhibit abrupt variations. The observed patterns follow what was observed in the two-dimensional case, where the polytopic
approach achieves lower iterations counts than AMG and Block-Jacobi for all the polynomial degrees. In particular, for $p\!=\!2$, the agglomerated multigrid
preconditioner again exhibits robust behavior, with iteration counts oscillating between $6$ and $11$ for most of the simulation. When $p\!=\!3$, iterations range from $8$ to $13$, far below the counts of AMG, which reach up to $21$ iterations, and of Block-Jacobi,
which reach up to $68$ iterations.

\begin{figure}[!t]
    \centering

    \includegraphics[width=\textwidth]{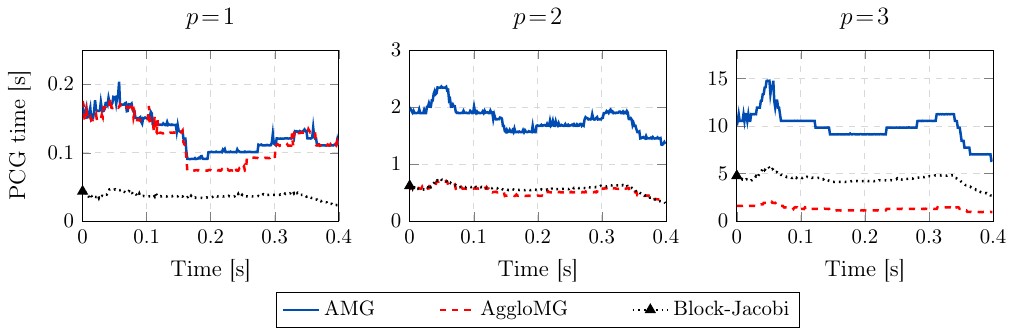}

    \caption{PCG iteration time per time step for Problem~\eqref{eqn:monodomain}, comparing AMG and agglomerated multigrid (AggloMG) for polynomial degrees $p\!=\!1,2,3$ for the three-dimensional
        ventricle test.}
    \label{fig:iteration_times_3D_monodomain}
\end{figure}

\subsection{Performance results on realistic three-dimensional test}
We now examine the performance of all the considered preconditioners for the realistic three-dimensional ventricle test case. To this end, we do not only consider iteration counts
as a metric, but also the wall-clock time spent in the PCG iterations, as well as the total time-to-solution of the whole simulation.
The simulations have been run on the MareNostrum 5 cluster at Barcelona Supercomputing Center, where each compute node features two Intel Sapphire Rapids 8460Y+ processors (80 cores per node at 2.3\si{\giga\hertz}) and 512\si{\giga\byte} of main memory. All nodes have
AVX-512 ISA extension, which allows 8 double-precision or 16 single-precision numbers to be processed per instruction when employing the matrix-free kernels detailed in Section~\ref{subsec:matrix-free}. The timings reported in the forthcoming figures and tables are in seconds.

Since our approach builds the coarse operators through the inherited strategy, we quantify the memory footprint of the proposed multigrid hierarchies and the cost of V-cycle application, using
as metrics the operator complexity, defined as:
\begin{equation}
    C_{\text{op}} \coloneqq \frac{\sum_{l=0}^{L-1} \text{nnz}(A_l)}{\text{nnz}(A_0)} \!>\! 1,
\end{equation}which is a customary measure considered in AMG methodologies~\cite{FalgoutIntroductionAMG}. In the definition above, $\text{nnz}(A_l)$ denotes the number of non-zero entries of the operator $A_l$ at level $l$.

Fig.~\ref{fig:iteration_times_3D_monodomain} displays the wall-clock time of each PCG iteration for all preconditioners when considering the whole simulation up to $T\!=\!0.4\si{\second}$. The plots highlight the higher efficiency of the agglomeration-based
multigrid preconditioner against the chosen AMG implementation. The agglomeration-based approach achieves consistently lower times for all polynomial degrees, yielding compatible performance also with linear elements, where
AMG is known to be very competitive. Despite the higher number of iterations, Block-Jacobi delivers the lowest total time-to-solution for linear elements, owing to its negligible per-iteration cost. With quadratic elements, the gap between
agglomerated multigrid and AMG is even more significant. Block-Jacobi yields times which are slightly higher than agglomeration-based multigrid. When cubic elements are considered, the performance gap between agglomerated multigrid and Block-Jacobi
is evident, with the latter showing significantly higher times, but better performance than AMG.

In Table~\ref{tab:timing_comparison}, we report some statistics about the PCG iterations for all preconditioners and polynomial degrees. From there, it can be appreciated
that with linear elements the time per PCG iteration is lower for the Block-Jacobi preconditioner, whereas AMG and agglomerated multigrid show higher times. However, the agglomeration-based multigrid preconditioner shows a speedup
with respect to AMG and Block-Jacobi that increases with the polynomial degree $p$: when $p\!=\!2$ it reaches a speedup factor greater than $3$ against AMG. When $p\!=\!3$, it achieves a factor of $3$ against Block-Jacobi and a factor of
roughly $7.5$ against AMG, yielding a significant improvement in the overall time-to-solution of the whole simulation.

\begin{figure}[!p]
    \centering

    \includegraphics[width=\textwidth]{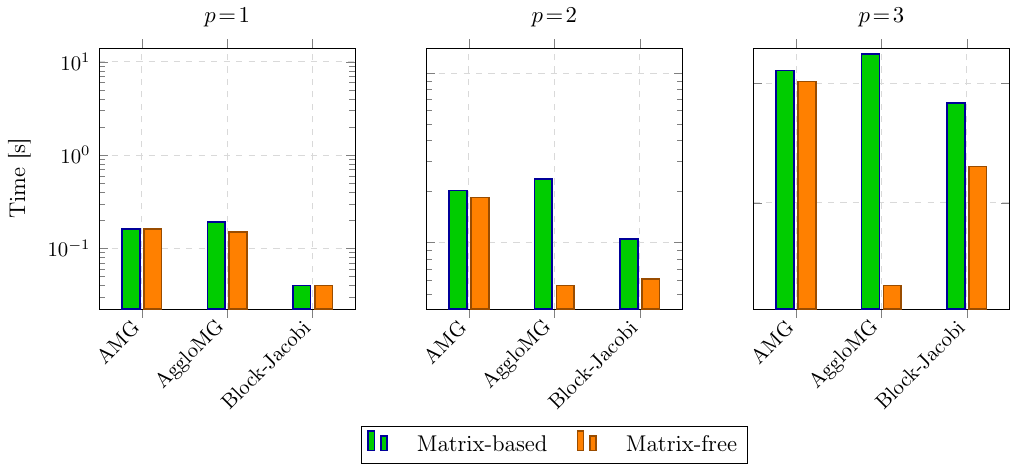}

    \caption{Comparison of average wall-clock time (measured in seconds) to perform one time step for polynomial degrees $p\!=\!1$, $p\!=\!2$, and $p\!=\!3$ on the realistic left ventricle mesh with $256$ MPI processes. Blue bars represent matrix-based variants, orange bars represent results obtained with matrix-free evaluations of the
        finest operator.}
    \label{fig:polynomial_comparison_MB_vs_MF}

    \vspace{1em}

    \centering

    \includegraphics[width=\textwidth]{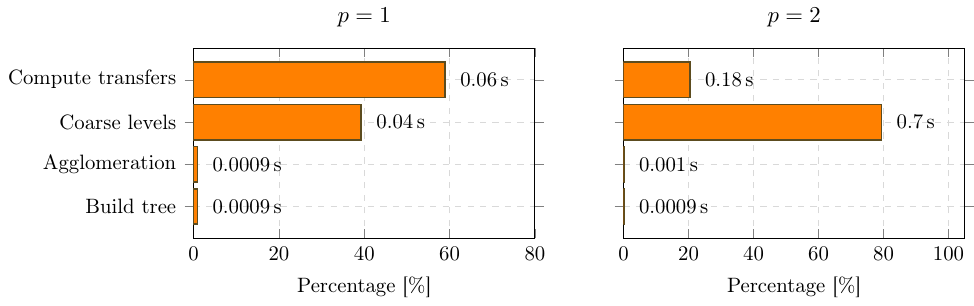}

    \par\vspace{0.3cm}
    \centering

    \begin{tikzpicture}
        \pgfmathsetmacro{\timeDataStructureC}{0.001}
        \pgfmathsetmacro{\timeCoarserLevelsC}{9.5}
        \pgfmathsetmacro{\timeAgglomeratesC}{0.001}
        \pgfmathsetmacro{\timeTransfersC}{0.64}
        \pgfmathsetmacro{\totalTimeC}{\timeDataStructureC + \timeCoarserLevelsC + \timeAgglomeratesC + \timeTransfersC}
        \pgfmathsetmacro{\percDataStructureC}{100 * \timeDataStructureC / \totalTimeC}
        \pgfmathsetmacro{\percCoarserLevelsC}{100 * \timeCoarserLevelsC / \totalTimeC}
        \pgfmathsetmacro{\percAgglomeratesC}{100 * \timeAgglomeratesC / \totalTimeC}
        \pgfmathsetmacro{\percTransfersC}{100 * \timeTransfersC / \totalTimeC}

        \begin{axis}[
                xbar,
                width=0.46\textwidth,
                height=0.3\textwidth,
                title={$p=3$},
                title style={font=\small},
                xlabel={Percentage [\%]},
                symbolic y coords={Build tree, Agglomeration, Coarse levels, Compute transfers},
                ytick=data,
                nodes near coords,
                nodes near coords align={horizontal},
                xmin=0, xmax=115,
                bar width=0.6cm,
                enlarge y limits=0.25,
                tick label style={font=\footnotesize},
                label style={font=\footnotesize},
                every node near coord/.append style={font=\footnotesize, anchor=west, xshift=3pt},
                point meta=explicit symbolic,
                nodes near coords style={/pgf/number format/.cd, fixed, precision=3, /tikz/.cd},
                grid=major,
                grid style={dashed, gray!30},
                clip=false
            ]
            \addplot[fill=orange, draw=teal!40!orange!60!black] coordinates {
                    (\percTransfersC,Compute transfers) [\timeTransfersC\,s]
                    (\percCoarserLevelsC,Coarse levels) [\timeCoarserLevelsC\,s]
                    (\percAgglomeratesC,Agglomeration) [\timeAgglomeratesC\,s]
                    (\percDataStructureC,Build tree) [\timeDataStructureC\,s]
                };
        \end{axis}
    \end{tikzpicture}
    \caption{Breakdown of the setup phase cost for the agglomeration-based multigrid preconditioner on the three-dimensional ventricle test case with $p\!=\!1$ (top left), $p\!=\!2$ (top right), and $p\!=\!3$ (bottom) with $256$ MPI processes. The percentages represent the relative time spent in each component.}
    \label{fig:setup_costs}
\end{figure}

The performance of the preconditioners when employing either matrix-based or matrix-free operator evaluations within the Krylov solver is compared in Fig.~\ref{fig:polynomial_comparison_MB_vs_MF}. Specifically, we show the time to perform a single PCG iteration
for polynomial degrees $p\!=\!1$, $p\!=\!2$, and $p\!=\!3$. In almost all instances, matrix-free evaluations (orange bars) improve the performance of the PCG iteration, compared to their matrix-based counterparts (green bars). When $p\!=\!1$, the speedup for agglomerated multigrid with matrix-free actions (instead of matrix-based ones) is
moderate, while essentially no gain is observed for AMG and Block-Jacobi with matrix-free operator evaluations. Matrix-based variants show higher times, which are instead in favor of Block-Jacobi.

The gain is substantially higher when $p\!=\!2$, with the agglomeration-based multigrid reaching a speedup factor greater than $4$ with respect to its matrix-based variant. The AMG and Block-Jacobi preconditioners with matrix-free operator evaluations show
a more pronounced speedup with respect to the matrix-based version. Such factors are, however, much smaller than the speedup
obtained for the agglomeration-based multigrid. When $p\!=\!3$, the agglomeration-based multigrid with matrix-free actions outperforms its matrix-based variant by a factor of $9$. It also shows a speedup greater than $3$ compared to Block-Jacobi. Another observation that can be drawn from Fig.~\ref{fig:polynomial_comparison_MB_vs_MF} is that
the performance gap between the agglomeration-based multigrid (with matrix-free evaluations) and AMG (with matrix-based evaluations, used in a classical black-box fashion) is quite significant, with a speedup factor of $3.6$ for $p\!=\!2$ and $7.8$ for $p\!=\!3$. Comparing in the same way to
simple Block-Jacobi with matrix-based operator evaluations, the speedup is $1.7$ for $p\!=\!2$ and $5.8$ for $p\!=\!3$.

Regarding the setup phase (column $t_{\text{setup}}$ in Table~\ref{tab:timing_comparison}), the chosen
AMG and Block-Jacobi implementations show a globally lower cost than the agglomeration-based one. Besides, this cost is amortized during the time stepping procedure, where the preconditioner is recycled at each time step.

To better understand the composition of the setup phase for the agglomeration-based preconditioner, we report in Fig.~\ref{fig:setup_costs} the breakdown of the setup cost into its main components: building the R-tree data structure from the fine mesh $\mathcal{T}_h$, traversing the tree hierarchy to create agglomerates, storing transfer
matrices $\{\mathcal{P}_{l+1}^l\}_{l\geq 0}$, and computing the coarser operators $\{A_l\}_{l=1}^{L-1}$ through the Galerkin product in Equation~\eqref{eqn:inherited_matrices}. It is evident from the
breakdown that the most expensive operation is the computation of the Galerkin triple product, which involves sparse matrix-matrix multiplications. This is particularly evident when increasing the polynomial degree to $p\!=\!2$ and $p\!=\!3$, with 80\% or more of the total setup time being spent in this phase. The second most expensive operation is the computation and storage of the intergrid transfer operators,
reaching almost 20\% of the total setup time for $p\!=\!2$ and 10\% for $p\!=\!3$. On the other hand, building the R-tree and creating agglomerates are totally inexpensive operations, taking together less than 1\% of the total setup time. Notably, this computational kernel depends \emph{only} on the geometry, i.e. the initial mesh $\mathcal{T}_h$, and not on the polynomial degree $p$.

\begin{figure}[!t]
    \centering

    \includegraphics[width=\textwidth]{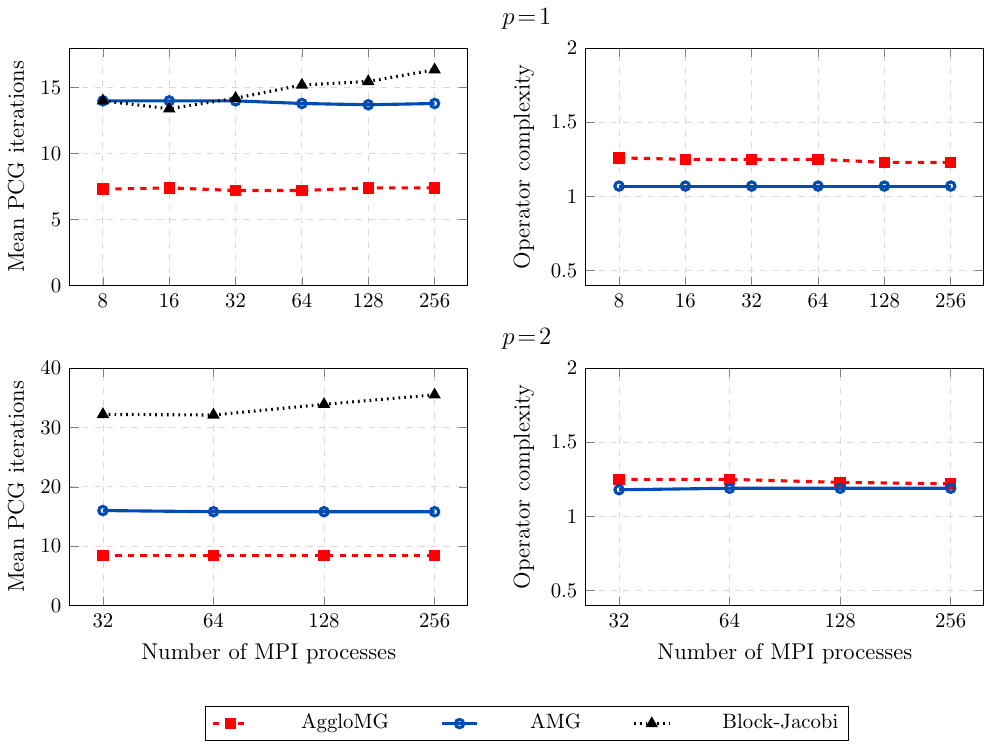}

    \caption{Scalability analysis with respect to the number of MPI processes for the three-dimensional ventricle test case with polynomial degrees $p\!=\!1$ (top row) and $p\!=\!2$ (bottom row). Left column: average number of PCG iterations per time step. The number of iterations is averaged over the first 20 time steps. Right column: operator complexity defined as the ratio of total non-zeros in the multigrid hierarchy to the non-zeros in the fine-grid operator. Notice that the operator complexity is a metric defined only for multigrid approaches, and hence it is not
        reported for Block-Jacobi.}
    \label{fig:mpi_iterations_operator_complexity}
\end{figure}

We detail in Fig.~\ref{fig:mpi_iterations_operator_complexity} the dependence of all preconditioners on the number of MPI processes employed. We report, for increasing numbers of processes, the average number
of PCG iterations\footnote{Averaged over the first twenty time steps.} per time step (left) and the operator complexity\footnote{The operator complexity for AMG has been computed using the interface exposed by \textsc{TrilinosML}.} (right). The first
row shows the results for $p\!=\!1$, while the second row for $p\!=\!2$. The number of PCG iterations for AMG and agglomeration-based multigrid is
clearly stable with respect to the number of processes involved, with the agglomeration-based multigrid showing lower iteration counts. Iteration counts for Block-Jacobi show an expected increase with the number of processes due to its limited algorithmic scalability. The computed
complexities for $p\!=\!1$ are essentially independent of the number of processes for both multilevel approaches, with AMG showing values very close to one, while the agglomeration-based multigrid
shows a constant value of $1.2$ and remains nicely bounded. Analogous behavior is observed for $p\!=\!2$: although iteration counts are higher than for $p\!=\!1$, both AMG
and agglomeration-based multigrid remain stable with respect to the number of processes, while Block-Jacobi shows again a slight increase. Operator complexities remain stable, though with a smaller gap; AMG shows slightly higher values for the operator complexity, around $1.2$. Notably,
the agglomeration-based multigrid is not affected by the variation in polynomial degree, showing
again an essentially constant operator complexity of $1.2$.

\begin{figure}[!t]
    \centering

    \begin{subfigure}{0.48\textwidth}
        \includegraphics[width=\textwidth]{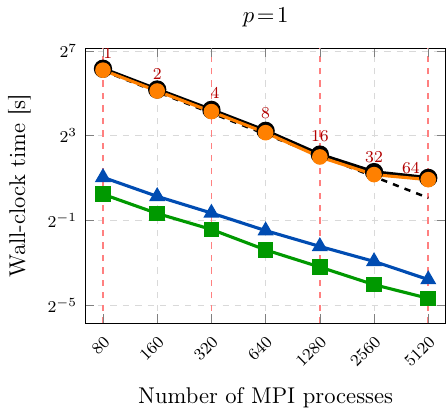}
    \end{subfigure}
    \hfill
    \begin{subfigure}{0.48\textwidth}
        \includegraphics[width=\textwidth]{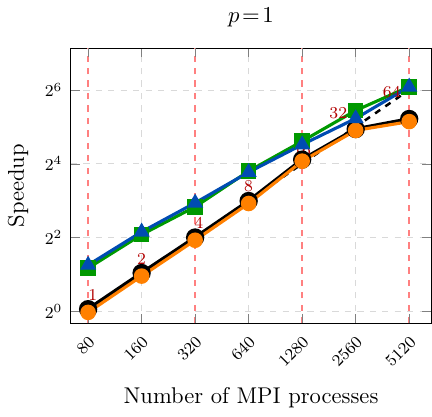}
    \end{subfigure}

    \vspace{0.4cm}
    \begin{subfigure}{0.48\textwidth}
        \includegraphics[width=\textwidth]{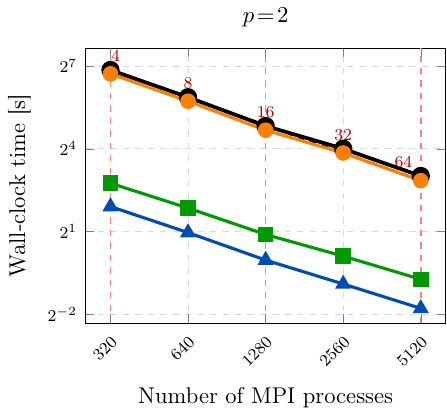}
    \end{subfigure}
    \hfill
    \begin{subfigure}{0.48\textwidth}
        \includegraphics[width=\textwidth]{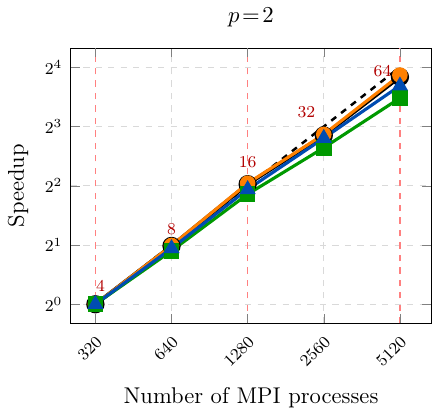}
    \end{subfigure}

    \vspace{0.2cm}
    \centering
    \begin{tikzpicture}
        \begin{axis}[
                hide axis,
                xmin=0, xmax=1, ymin=0, ymax=1,
                legend columns=-1,
                legend style={column sep=0.4cm, text=black, font=\footnotesize},
            ]
            \addlegendimage{black, dashed, line width=1.2pt}
            \addlegendentry{Ideal}
            \addlegendimage{black, solid, line width=1.8pt, mark=*, mark size=3.5pt, mark options={solid, fill=black}}
            \addlegendentry{Total}
            \addlegendimage{orange, solid, line width=1.5pt, mark=*, mark size=3pt, mark options={solid, fill=orange}}
            \addlegendentry{Linear solver}
            \addlegendimage{green!60!black, solid, line width=1.5pt, mark=square*, mark size=3pt, mark options={solid, fill=green!60!black}}
            \addlegendentry{Assembly}
            \addlegendimage{teal!60!blue, solid, line width=1.5pt, mark=triangle*, mark size=3pt, mark options={solid, fill=teal!60!blue}}
            \addlegendentry{Rhs assembly}
        \end{axis}
    \end{tikzpicture}
    \caption{Strong scaling results for the realistic ventricle test case with $p\!=\!1$ (top) and $p\!=\!2$ (bottom), using 23,937,408 and 80,788,752 degrees of freedom, respectively. Left: wall-clock time for different simulation components. Right: speedup computed relative to the baseline configuration with 80 processes for $p\!=\!1$ and 320 processes for $p\!=\!2$. Node counts are annotated in red.}
    \label{fig:strong_scaling}
\end{figure}

Finally, we perform a strong scaling test for the agglomeration-based multigrid preconditioner. We consider the realistic left ventricle test case with polynomial degree $p\!=\!1$, using a mesh $\mathcal{T}_h$ comprising 2,992,176 elements, which corresponds to 23,937,408 degrees of freedom. We run
the simulation for twenty time steps with $\Delta t\!=\!10^{-4}\si{\second}$, varying the number of
MPI processes from $80$ to 5,120 (i.e., from $1$ to $64$ compute nodes). Fig.~\ref{fig:strong_scaling} (top row) shows the strong scaling behavior of several
components of the simulations. In particular, the linear solver exhibits ideal scaling up to $32$ nodes, with efficiency remaining above $90\%$ up to $32$ nodes (2,560 processes), demonstrating the good parallel scalability of the proposed approach. Notice that at the highest number of processes ($64$ nodes, $5,120$ processes), less
than $10,000$ degrees of freedom are assigned to each process, hence justifying the mild degradation due to communication overhead, which
matches previous observations in the literature~\cite{AlgorithmsDealII}. The other components of the simulation, such as assembly of the system matrix (performed only once), and right-hand side assembly (performed at every time-step), also show ideal scaling behavior. In the bottom row of Fig.~\ref{fig:strong_scaling}, we repeat the
experiment with $p\!=\!2$, corresponding to 80,788,752 degrees of freedom, from $4$ up to $64$ compute nodes. The linear solver shows almost ideal scaling up to $64$ nodes, confirming the good scalability properties of the
agglomeration-based multigrid preconditioner also for higher polynomial degrees.

All in all, we believe the series of experiments above shows the potential of the agglomeration-based preconditioner, particularly when comparing against highly optimized
and widely used solvers. The high efficiency of our aggregation strategy suggests that it could be readily integrated into existing high-performance AMG frameworks (such as \textsc{psctoolkit}~\cite{DAMBRA2023100463}), in order
to provide \emph{geometrically informed} AMG hierarchies.
Nevertheless, our implementation can be optimized further: within the Galerkin triple product computation, and during the evaluation
of coarser operators. We leave these improvements to future work. Despite these possible enhancements, we conclude
this discussion by observing that the agglomeration-based preconditioner is an efficient and competitive solver both in
terms of iteration counts and time-to-solution for the present problem, with its advantage over AMG and Block-Jacobi growing with the polynomial degree.

\begin{table}[!h]
    \centering
    \begin{tabular}{@{}lcccccc@{}}
        \toprule
        $p$ & Preconditioner & $t_{\text{setup}}$ [s] & It. & $t_{\text{CG}}$ [s] & $t_{\text{TOT}}$ [s] & DoFs                        \\
        \midrule
        \multirow{3}{*}{$1$}
            & AMG            & 0.11                   & 13  & 0.12                & $5.68\times 10^2$    & \multirow{3}{*}{2,992,176}  \\
            & AggloMG        & 0.27                   & 6   & 0.11                & $4.83\times 10^2$    &                             \\
            & Block-Jacobi   & 0.05                   & 16  & 0.03                & $1.60\times 10^2$    &                             \\
        \midrule
        \multirow{3}{*}{$2$}
            & AMG            & 0.68                   & 15  & 1.78                & $7.19\times 10^3$    & \multirow{3}{*}{10,098,594} \\
            & AggloMG        & 1.07                   & 8   & 0.52                & $2.14\times 10^3$    &                             \\
            & Block-Jacobi   & 0.70                   & 34  & 0.57                & $2.33\times 10^3$    &                             \\
        \midrule
        \multirow{3}{*}{$3$}
            & AMG            & 4.36                   & 14  & 10.0                & $4.01\times 10^4$    & \multirow{3}{*}{23,937,408} \\
            & AggloMG        & 10.8                   & 8   & 1.32                & $5.42\times10^3$     &                             \\
            & Block-Jacobi   & 5.9                    & 52  & 4.39                & $1.76\times 10^4$    &                             \\
        \bottomrule
    \end{tabular}
    \caption{Timing comparison between AMG, agglomerated multigrid (AggloMG), and Block-Jacobi preconditioners for polynomial degrees $p=1,2,3$ on the three-dimensional ventricle test with 256 MPI processes. Wall-clock times for the setup phase ($t_{\text{setup}}$), the average wall-clock time per time step ($t_{\text{CG}}$), and the overall simulation time ($t_{\text{TOT}}$) are reported, along with the average number of PCG iterations per time step (It.), rounded to the nearest integer. The total number of time steps is 4000.}

    \label{tab:timing_comparison}
\end{table}

\begin{table}[!t]
    \centering
    \begin{minipage}{0.48\textwidth}
        \centering
        \footnotesize
        \begin{tabular}{@{}lcrc@{}}
            \toprule
            \multicolumn{4}{c}{Hierarchy for ventricle mesh ($128$ processes)}                              \\
            \midrule
            Level $l$ & Card($\mathcal{T}_l$) & DoFs      & Card($\mathcal{T}_{l-1}$)/Card($\mathcal{T}_l$) \\
            \midrule
            $l=0$     & 374,022               & 2,992,176 &                                                 \\
            $l=1$     & 46,809                & 374,472   & $\sim 7.9$                                      \\
            $l=2$     & 5,908                 & 47,264    & $\sim 7.9$                                      \\
            $l=3$     & 768                   & 6,144     & $\sim 7.7$                                      \\
            \bottomrule
        \end{tabular}
    \end{minipage}%
    \hspace{0.04\textwidth}%
    \begin{minipage}{0.48\textwidth}
        \centering
        \footnotesize
        \begin{tabular}{@{}crc@{}}
            \toprule
            \multicolumn{3}{c}{Hierarchy for ventricle mesh ($256$ processes)}                  \\
            \midrule
            Card($\mathcal{T}_l$) & DoFs      & Card($\mathcal{T}_{l-1}$)/Card($\mathcal{T}_l$) \\
            \midrule
            374,022               & 2,992,176 &                                                 \\
            46,871                & 374,968   & $\sim 7.9$                                      \\
            5,987                 & 47,896    & $\sim 7.8$                                      \\
            768                   & 6,144     & $\sim 7.8$                                      \\
            \bottomrule
        \end{tabular}
    \end{minipage}
    \caption{Coarsened hierarchies $\{\mathcal{T}_l\}_{l=0}^3$ on top of the original mesh $\mathcal{T}_0\equiv \mathcal{T}_h$, partitioning the realistic ventricle mesh across different numbers of processes. The number of agglomerates and degrees of freedom (with $p=1$) per level are shown in the second and third columns, while the ratio between the cardinality of consecutive grids is reported in the last column.}
    \label{tab:Ratio_meshes_ventricle}
\end{table}

\section{Conclusions and outlook}\label{sec:conc}
\noindent This work addressed the development of a novel multilevel preconditioner for the discontinuous Galerkin discretization of the monodomain problem in cardiac electrophysiology. The preconditioner exploits the flexibility of
discontinuous Galerkin methods in terms of their ability to handle general polytopic shapes, which can be generated by agglomeration of fine mesh elements. The R-tree-based agglomeration algorithm from~\cite{FEDER2025113773} enables the construction of nested hierarchies of coarser polytopic grids
and their associated operators, yielding a robust and computationally efficient multigrid preconditioner. Our preconditioner exhibits good scalability with respect to both the
polynomial degree and the number of levels, maintaining low iteration counts for a cohort of simulations. We employ state-of-the-art matrix-free operator evaluation to further improve computational efficiency. The effectiveness of the proposed methodology and its performance are validated through a sequence of tests involving two-dimensional and three-dimensional domains, different ionic models, and employing real unstructured geometries. Future work will focus on
efficient coarse-operator evaluation, which is currently under investigation. Integration with existing high-performance
AMG packages and application to more realistic settings and problems, including agglomeration of adaptively refined meshes, are also topics
for future research.

\section*{Data availability}
The source code implementing the methodology and the experiments showcased in this paper are publicly available at the GitHub page of \textsc{polyDEAL}~\cite{polydeal}. Detailed instructions for compilation and execution are provided in the repository.

\section*{CRediT author statement}
\textbf{M. Feder}: Formal analysis, Investigation, Methodology, Software, Visualization, Writing -- Original draft, Writing -- review and editing. \textbf{P.C. Africa}: Conceptualization, Formal analysis, Methodology, Supervision, Writing -- review and editing.

\section*{Acknowledgments}
\noindent MF and PCA acknowledge support of the European Research Council (ERC) under the European Union's Horizon 2020 research and innovation programme (call HORIZON-EUROHPC-JU-2023-COE-03, grant agreement No. 101172493 "dealii-X"). PCA acknowledges support by the
European Union - NextGenerationEU, in the framework of the iNEST - Interconnected Nord-Est Innovation Ecosystem (iNEST ECS00000043 - CUP G93C22000610007) and its CC5 Young
Researchers initiative. The authors thank the Barcelona Supercomputing Center for the computational resources made available through the EuroHPC project EHPC-DEV-2026D02-275. The authors are members of Gruppo Nazionale per il Calcolo Scientifico (GNCS), Istituto Nazionale di Alta Matematica (INdAM). The authors are thankful to Prof. Andrea Cangiani, Prof. Luca Heltai, and Prof. Fabio Durastante for useful discussions.

\appendix
\section{Multigrid and model parameters}
\label{app:parameters}
\noindent This appendix reports the parameters employed in the multigrid preconditioners and ionic models used in the numerical tests.

\begin{table}[h]
    \centering
    \begin{tabular}{ll}
        \toprule
        Parameter             & Value                \\
        \midrule
        Smoother              & Chebyshev            \\
        Smoother degree       & 3                    \\
        Smoother sweeps       & 3                    \\
        V-cycle applications  & 1                    \\
        Aggregation type      & Smoothed aggregation \\
        Aggregation threshold & \(0.2\)              \\
        Max size coarse level & 2000                 \\
        Coarse solver         & Amesos-KLU           \\
        \bottomrule
    \end{tabular}
    \caption{Parameters for TrilinosML preconditioner~\cite{Trilinos} (Trilinos 14.4.0). When using polynomial degrees higher than one, $\mathtt{higher\_order\_elements}$ is set to true. In that case, the $\texttt{Uncoupled}$ option is automatically
        selected as aggregation type.}
    \label{tab:amg_params}
\end{table}

\begin{table}[H]
    \centering
    \begin{minipage}{0.48\textwidth}
        \centering
        \begin{tabular}{@{}lc@{}}
            \toprule
            \textbf{Parameter}   & \textbf{Value}       \\
            \midrule
            Cycle type           & V-cycle              \\
            Number of cycles     & 1                    \\
            Smoother             & Chebyshev (degree 3) \\
            Pre-smoothing steps  & 3                    \\
            Post-smoothing steps & 3                    \\
            Coarse solver        & \textsc{MUMPS}       \\
            \bottomrule
        \end{tabular}
        \caption{Agglomerated multigrid configuration for the two-dimensional test case.}
        \label{tab:agglo_mg_config_2d}
    \end{minipage}
    \hfill
    \begin{minipage}{0.48\textwidth}
        \centering
        \begin{tabular}{@{}lc@{}}
            \toprule
            \textbf{Parameter}   & \textbf{Value}       \\
            \midrule
            Cycle type           & V-cycle              \\
            Number of cycles     & 1                    \\
            Smoother             & Chebyshev (degree 3) \\
            Pre-smoothing steps  & 3                    \\
            Post-smoothing steps & 3                    \\
            Coarse solver        & PCG+AMG              \\
            \bottomrule
        \end{tabular}
        \caption{Agglomerated multigrid configuration for the three-dimensional test case.}
        \label{tab:agglo_mg_config_3d}
    \end{minipage}
\end{table}

\begin{table}[h]
    \centering
    \begin{tabular}{@{}lr@{}}
        \toprule
        Parameter    & \multicolumn{1}{c}{Value}                                      \\
        \midrule
        $C_m$        & $1\times10^{-2}$ [\si{\meter^{-1}}]                            \\
        $\chi_m$     & $1\times10^{5}$ [\si{\farad\per\meter\squared}]                \\
        $\Sigma$     & $1.2\times10^{-1}$ [\si{\meter\squared\per\second\squared}]    \\
        $\mathbb{D}$ & $\Sigma I_{2\times2}$ [\si{\meter\squared\per\second\squared}] \\
        $\kappa$     & $1.95\times10^{1}$                                             \\
        $\epsilon$   & $4\times10^{1}$                                                \\
        $\Gamma$     & $1\times10^{-1}$                                               \\
        $a$          & $1.3\times10^{-2}$                                             \\
        \bottomrule
    \end{tabular}
    \caption{FitzHugh-Nagumo parameters used in the 2D numerical test.}
    \label{tab:FHZ_params}
\end{table}

\begin{table}[h]
    \centering
    \small
    \begin{tabular}{@{}lrlrlr@{}}
        \toprule
        Parameter        & Value                                                       & Parameter           & Value                                     & Parameter            & Value                \\
        \midrule
        $\sigma_l$       & $1.0\times10^{-4}$ [\si{\meter\squared\per\second\squared}] & $\tau_{1}^{'}$      & $6\times10^{-2}$ [\si{\second^{-1}}]      & $k_{2}$              & $6.5\times10^{1}$    \\
        $\sigma_t$       & $0.5\times10^{-4}$ [\si{\meter\squared\per\second\squared}] & $\tau_{1}^{''}$     & $1.15$ [\si{\second^{-1}}]                & $k_{3}$              & $2.0994$             \\
        $\sigma_n$       & $0.1\times10^{-4}$ [\si{\meter\squared\per\second\squared}] & $\tau_{2}^{'}$      & $7\times10^{-2}$ [\si{\second^{-1}}]      & $k_{so}$             & $2.0$                \\
        $\tau_{o}^{'}$   & $6\times10^{-3}$ [\si{\second^{-1}}]                        & $\tau_{2}^{''}$     & $2\times10^{-2}$ [\si{\second^{-1}}]      & $w_{\infty}^{\star}$ & $9.4\times10^{-1}$   \\
        $\tau_{o}^{''}$  & $6\times10^{-3}$ [\si{\second^{-1}}]                        & $\tau_{3}^{'}$      & $2.7342\times10^{-3}$ [\si{\second^{-1}}] & $V_1$                & $3\times10^{-1}$     \\
        $\tau_{so}^{'}$  & $4.3\times10^{-2}$ [\si{\second^{-1}}]                      & $\tau_{2}^{''}$     & $2\times10^{-2}$ [\si{\second^{-1}}]      & $w_{\infty}^{\star}$ & $9.4\times10^{-1}$   \\
        $\tau_{so}^{''}$ & $2\times10^{-4}$ [\si{\second^{-1}}]                        & $\tau_{3}^{'}$      & $2.7342\times10^{-3}$ [\si{\second^{-1}}] & $V_1$                & $3\times10^{-1}$     \\
        $\tau_{si}$      & $2.8723\times10^{-3}$ [\si{\second^{-1}}]                   & $\tau_{3}^{''}$     & $3\times10^{-3}$                          & $V_{1m}$             & $1.5\times10^{-2}$   \\
        $\tau_{fi}$      & $1.1\times10^{-4}$ [\si{\second^{-1}}]                      & $\tau_{2}^{+}$      & $2.8\times10^{-1}$ [\si{\second^{-1}}]    & $V_2$                & $1.5\times10^{-2}$   \\
        $\tau_{1}^{+}$   & $1.4506\times10^{-3}$ [\si{\second^{-1}}]                   & $\tau_{2}^{\infty}$ & $7\times10^{-2}$ [\si{\second^{-1}}]      & $V_{2m}$             & $3\times10^{-2}$     \\
                         &                                                             &                     &                                           & $V_3$                & $9.087\times10^{-1}$ \\
                         &                                                             &                     &                                           & $\hat{V}$            & $1.58$               \\
                         &                                                             &                     &                                           & $V_{o}$              & $6\times10^{-3}$     \\
                         &                                                             &                     &                                           & $V_{so}$             & $6.5\times10^{-1}$   \\
        \bottomrule
    \end{tabular}
    \caption{Parameters for Bueno-Orovio model used in the 3D numerical test with the idealized ventricle mesh. In this case $\chi_m\equiv C_m =1$. In the realistic ventricle test, the parameters are the same except for the conductivities, which are set to $\sigma_l=10\times10^{-5}$ [\si{\meter\squared\per\second\squared}], $\sigma_t=2\times10^{-5}$ [\si{\meter\squared\per\second\squared}], and $\sigma_n=2\times10^{-5}$ [\si{\meter\squared\per\second\squared}]. Values have been chosen according to~\cite{lifex-ep}.}
    \label{tab:BO_params}
\end{table}

\begin{table}[h!]
    \centering
    \resizebox{\textwidth}{!}{%
        \begin{tabular}{@{}ccc cccc cccc cccc@{}}
            \toprule
            \multirow{2}{*}{$p$} & \multirow{2}{*}{DoFs}                                      & \multirow{2}{*}{$\Delta t$ [s]}
                                 & \multicolumn{4}{c}{PCG Iterations (\textbf{AggloMG})}
                                 & \multicolumn{4}{c}{PCG Iterations (\textbf{AMG})}
                                 & \multicolumn{4}{c}{PCG Iterations (\textbf{Block-Jacobi})}                                                                                                              \\
            \cmidrule(lr){4-7} \cmidrule(lr){8-11} \cmidrule(lr){12-15}
                                 &                                                            &
                                 & Mean                                                       & Std Dev                         & Min   & Max
                                 & Mean                                                       & Std Dev                         & Min   & Max
                                 & Mean                                                       & Std Dev                         & Min   & Max                                                              \\
            \midrule
            \multirow{3}{*}{$1$}
                                 & \multirow{3}{*}{3,263,232}                                 & $5\times 10^{-4}$               & 10.48 & 0.51 & 9  & 11 & 19.89 & 0.98 & 17 & 22 & 24.26 & 2.18 & 21 & 28 \\
                                 &                                                            & $1\times 10^{-4}$               & 7.78  & 0.58 & 6  & 9  & 14.70 & 0.94 & 11 & 17 & 13.53 & 1.27 & 11 & 16 \\
                                 &                                                            & $5\times 10^{-5}$               & 6.98  & 0.38 & 6  & 8  & 13.00 & 0.62 & 10 & 15 & 10.43 & 0.97 & 9  & 12 \\
            \midrule
            \multirow{3}{*}{$2$}
                                 & \multirow{3}{*}{11,013,408}                                & $5\times 10^{-4}$               & 10.90 & 0.69 & 9  & 12 & 19.92 & 1.21 & 17 & 22 & 49.91 & 3.70 & 43 & 56 \\
                                 &                                                            & $1\times 10^{-4}$               & 8.62  & 0.63 & 7  & 10 & 15.94 & 0.89 & 14 & 18 & 28.95 & 2.33 & 24 & 33 \\
                                 &                                                            & $5\times 10^{-5}$               & 7.39  & 0.59 & 5  & 9  & 14.07 & 0.91 & 12 & 16 & 21.89 & 1.97 & 18 & 25 \\
            \midrule
            \multirow{3}{*}{$3$}
                                 & \multirow{3}{*}{26,105,856}                                & $5\times 10^{-4}$               & 12.21 & 1.11 & 11 & 14 & 17.48 & 1.23 & 15 & 19 & 71.86 & 4.88 & 63 & 81 \\
                                 &                                                            & $1\times 10^{-4}$               & 8.89  & 0.67 & 7  & 10 & 12.10 & 0.89 & 10 & 14 & 44.71 & 3.53 & 38 & 50 \\
                                 &                                                            & $5\times 10^{-5}$               & 7.88  & 0.61 & 6  & 9  & 10.48 & 0.84 & 8  & 12 & 35.04 & 2.95 & 29 & 40 \\
            \bottomrule
        \end{tabular}}
    \caption{Statistics (mean value, standard deviation, and min-max) of PCG iterations per time step for the idealized left ventricle test (up to $T=0.1$ s) case across different time steps $\Delta t$ and polynomial degrees $p$, for the three preconditioners AggloMG, AMG, and Block-Jacobi.}
    \label{tab:multigrid_iterations_delta_t}
\end{table}

\bibliographystyle{abbrv}
\bibliography{refs}
\end{document}